\documentclass[11pt,draft]{article}
\usepackage{ourstyle-private}
\usepackage{snapshot}

\DeclareGraphicsExtensions{.pdf,.png}

\title{%
Benchmarking the Immersed Finite Element Method for Fluid-Structure Interaction Problems
}

\author{%
Saswati Roy\\
Department of Engineering Science and Mechanics\\
The Pennsylvania State University\\
212 Earth and Engineering Sciences Building\\
University Park
PA 16802
USA
\and
Luca Heltai\footnote{Corresponding Author. Email: Luca Heltai \texttt{<luca.heltai@sissa.it>}; Tel.: +39 040 3787-449}
\\
Scuola Internazionale Superiore di Studi Avanzati\\
                    Via Bonomea 265\\
                    34136 Trieste, Italy
\and
Francesco Costanzo\\
Center for Neural Engineering\\
The Pennsylvania State University\\
W-315 Millennium Science Complex\\
University Park
PA 16802
USA}

\begin{document}

\maketitle
               
\begin{abstract}
We present an implementation of a fully variational formulation of an
immersed method for fluid-structure interaction problems based on the
finite element method.  While typical implementation of immersed
methods are characterized by the use of approximate Dirac delta
distributions, fully variational formulations of the method do not
require the use of said distributions.  In our implementation the
immersed solid is general in the sense that it is not required to have
the same mass density and the same viscous response as the surrounding
fluid.  We assume that the immersed solid can be either viscoelastic
of differential type or hyperelastic.  Here we focus on the validation
of the method via various benchmarks for fluid-structure interaction
numerical schemes.  This is the first time that the interaction of
purely elastic compressible solids and an incompressible fluid is
approached via an immersed method allowing a direct comparison with
established benchmarks.
\end{abstract}


\allowdisplaybreaks{                             %


\section{Introduction}
Immersed methods for fluid structure interaction (\acro{FSI}) problems were pioneered by Peskin and his co-workers \citep{Peskin_1977_Numerical_0,Peskin_2002_The-immersed_0}.  They proposed an approach called the immersed boundary method (\acro{IBM}), in which the equations governing the fluid motion have body force terms describing the \acro{FSI}.  The equations are integrated via a finite difference (\acro{FD}) method and the body force terms are computed by modeling the solid body as a network of elastic fibers.  As such, this system of forces has singular support (the \emph{boundary} in the method's name) and is implemented via Dirac-$\delta$ distributions. The configuration of the fiber network is represented via a discrete set of points whose motion is then related to that of the fluid again via Dirac-$\delta$ distributions.  We should clarify that the fiber network in question can be configured so as to represent both thin elastic interfaces as well as thick elastic bodies.  In the numerical implementation of this method the Dirac-$\delta$ distributions are aproximated as \emph{functions}.  A recent paper by  \cite{Fai2014Immersed-Bounda-0} offers a detailed stability analysis in the context of problems with variable density and viscosity.

Immersed methods based on the finite element method (\acro{FEM}) has been formulated by various authors \citep{BoffiGastaldi_2003_A-Finite_0,WangLiu_2004_Extended_0,ZhangGerstenberger_2004_Immersed_0,BoffiGastaldiHeltaiPeskin-2008-a}.  \cite{BoffiGastaldi_2003_A-Finite_0} were the first to show that a variational approach to immersed methods does not necessitate the approximation of Dirac-$\delta$ distributions as they naturally disappears in the weak formulation.  The thrust of the work by \cite{WangLiu_2004_Extended_0} and \cite{ZhangGerstenberger_2004_Immersed_0} was to remove the requirement that the immersed solid be a fiber network.  They also included the ability to accommodate density differences between solid and fluid, as well as compressible materials in addition to incompressible ones. While they proposed an approach applicable to solid bodies of general topological and constitutive characteristics they maintained the use of approximated Dirac-$\delta$ distribution through a strategy called the reproducing kernel particle method (\acro{RKPM}).

Recently, \cite{Heltai2012Variational-Imp0} proposed a generalization of the approach by \cite{BoffiGastaldiHeltaiPeskin-2008-a} in which a fully variational \acro{FEM} formulation is shown to be applicable to problems with immersed bodies of general topological and constitutive characteristics and without the use of Dirac-$\delta$ distributions.  The discussion in \cite{Heltai2012Variational-Imp0} focused on the construction of natural interpolation operators between the fluid and the solid discrete spaces that guarantee semi-discrete stability estimates and strong consistency.  Since the formulation in \cite{Heltai2012Variational-Imp0} is applicable to solid bodies with pure hyperelastic behavior, i.e., without a viscous component in the stress response, in this paper we show that the method in question satisfies the benchmark tests by \cite{Turek2006Proposal-for-Nu0}.  This is an important result given that these benchmarks have become a \emph{de facto} standard in the FSI computational community, and given that previous immersed methods could not satisfy them due to intrinsic model restrictions.  In this sense, and to the best of our knowledge, our results are the first of their kind.  Along with these important results, we also illustrate the application of our method to a three-dimensional problem whose geometry is similar to that in the two-dimensional benchmark tests.

\section{Formulation}
\label{sec: Problem Formulation}

\subsection{Basic notation and governing equations}
\label{subsec: Basic notation and governing equations}
Referring to Fig.~\ref{fig: current_configuration}, 
\begin{figure}[htb]
    \centering
    \includegraphics{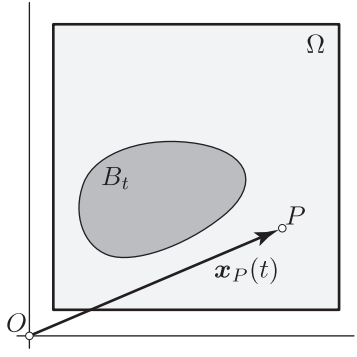}
    \caption{Current configuration $B_{t}$ of a body $\mathscr{B}$ immersed in a fluid occupying the domain $\Omega$.}
    \label{fig: current_configuration}
\end{figure}
$B_{t}$ is a body immersed in a fluid, the latter occupying $\Omega\setminus B_{t}$, where $\Omega$ is a fixed control volume.  The body's motion is described by a diffeomorphism $\bv{\zeta}: B \to B_{t}$, $\bv{x} = \bv{\zeta}(\bv{s},t)$, where $B$ is the body's reference configuration, $\bv{s} \in B$, $\bv{x} \in \Omega$, and time $t \in [0,T)$, with $T > 0$.  Away from boundaries, the motion $B_{t}$ and the fluid are both governed by the balance of mass and momentum, respectively,
\begin{gather}
\label{eq: Balance of mass and momentum}
\frac{\partial\rho}{\partial t} + \nabla \cdot (\rho \bv{u}) = 0
\quad \text{and} \quad 
\nabla \cdot \bv{\sigma} + \rho \bv{b} = \rho \biggl[\frac{\partial\bv{u}}{\partial t} + (\nabla \bv{u}) \bv{u} \biggr],
\end{gather}
where $\rho(\bv{x}, t)$ is the mass density, $\bv{u}(\bv{x},t)$ is the velocity, $\bv{\sigma}(\bv{x},t)$ is the Cauchy stress, $\bv{b}(\bv{x},t)$ is the external force density per unit mass, and where $\nabla$ and $(\nabla \cdot{})$ denote the gradient and divergence operators, respectively.  Equations~\eqref{eq: Balance of mass and momentum} hold both for the solid and the fluid, which can be distinguished via their constitutive equations. We assume that $\bv{u}(x,t)$ is continuous across $\partial B_{t}$, the boundary of $B_{t}$.  Along with the jump conditions for the momentum balance laws, this implies that the traction field is also continuous across $\partial B_{t}$. Equations~\eqref{eq: Balance of mass and momentum} are complemented by the following boundary conditions:
\begin{equation}
\label{eq: boundary conditions}
\bv{u}(\bv{x},t) = \bv{u}_{g}(\bv{x},t),\quad \text{for $\bv{x} \in \partial\Omega_{D}$,}
\quad \text{and} \quad
\bv{\sigma}(\bv{x},t) \bv{n}(\bv{x},t) = \mathfrak{s}_{g}(\bv{x},t), \quad \text{for $\bv{x} \in \partial\Omega_{N}$,}
\end{equation}
where $\bv{u}_{g}$ and $\mathfrak{s}_{g}$ are prescribed velocity and surface traction fields, $\bv{n}$ is the outward unit normal to $\partial \Omega$, and where $\partial\Omega_{D} \cup \partial\Omega_{N} = \partial\Omega$ and $\partial\Omega_{D} \cap \partial\Omega_{N} = \emptyset$.

\paragraph{Fluid's constitutive response.}
The fluid is assumed to be Newtonian with constant mass density $\rho_{\f}$ and stress
\begin{equation}
\label{eq: incompressible NS fluid}
\bv{\sigma} = -p \tensor{I} + \bv{\sigma}^{v}_{\f}
\quad \text{and} \quad 
\bv{\sigma}^{v}_{\f} = \mu_{\f} \bigl(\nabla \bv{u} + \trans{\nabla\bv{u}}\bigr),
\end{equation}
where $\f$ stands for `fluid,' $p$ is the pressure, $\tensor{I}$ is the identity tensor, $\bv{\sigma}^{v}_{\f}$ the linear viscous component of the stress, and $\mu_{\f} > 0$ is the dynamic viscosity.  For constant $\rho_{\f}$, the first of Eqs.~\eqref{eq: Balance of mass and momentum} reduces to $\nabla \cdot \bv{u} = 0$ ($\bv{x} \in \Omega\setminus B_{t}$) and $p$ is a multiplier for the enforcement of this constraint.

\paragraph{Solid's constitutive response.}
We consider both incompressible and compressible materials.  For the incompressible case, Cauchy stress is assumed to be
\begin{equation}
\label{eq: Cauchy Response Function}
\bv{\sigma} = -p \tensor{I} + \bv{\sigma}^{e}_{\s} + \bv{\sigma}^{v}_{\s},
\end{equation}
where $\s$ stands for `solid,' $p$ is a multiplier enforcing incompressibility, and where $\bv{\sigma}^{e}_{\s}$ and $\bv{\sigma}^{v}_{\s}$ are the elastic and viscous parts of the stress, respectively:
\begin{equation}
\label{eq: elastic and viscous part solid}
\bv{\sigma}^{e}_{\s} = J^{-1} \tensor{P}^{e}_{\s}\trans{\tensor{F}}, \quad
\tensor{P}^{e}_{\s} = \frac{\partial W^{e}_{\s}(\tensor{F})}{\partial \tensor{F}},
\quad \text{and} \quad 
\bv{\sigma}^{v}_{\s} =  \mu_{\s} \bigl(\nabla \bv{u} + \trans{\nabla\bv{u}}\bigr),
\end{equation}
where $W^{e}_{\s}$ is the elastic strain energy per unit referential volume, $\tensor{F} = \partial \bv{\zeta}(\bv{s},t)/\partial \bv{s}$ is the deformation gradient, $J = \det\tensor{F}$, $\tensor{P}^{e}_{\s}$ is the first Piola-Kirchhoff stress tensor, and $\mu_{\s}$ is the solid's dynamic viscosity.  The constitutive response for the compressible case is identical to that in Eq.~\eqref{eq: Cauchy Response Function} except for the multiplier $p$, which is not needed in this case.  Whether compressible or incompressible, we assume that $\mu_{\s} \geq 0$, that is, we assume that $\mu_{\s}$ might be equal to zero, in which case the solid is purely elastic.  We assume that $W_{\s}^{e}(\tensor{F})$ is a $C^{1}$ convex function over the set of second order tensor with positive determinant.  Finally, $\rho_{\s_{0}}(\bv{s})$ is the referential solid's mass density and we recall that the first of Eq.~\eqref{eq: Balance of mass and momentum} is equivalently expressed as $\rho_{\s_{0}}(\bv{s}) = \rho_{\s}(\bv{x},t)\big|_{\bv{x} = \bv{\zeta}(\bv{s},t)} J(\bv{s},t)$ ($\bv{s} \in B$).

\section{Variational and Discrete Formulations}
\label{sec: Reformulation of the governing equations}
The results in this paper have been obtained using the method in \cite{Heltai2012Variational-Imp0}.  For the sake of completeness, here we summarize the essential elements of this method.

Our formulation has a single velocity field representing the velocity of the fluid for $\bv{x} \in \Omega\setminus B_{t}$ and the velocity of the solid for $\bv{x} \in B_{t}$.  The motion of the solid is described by the displacement field $\bv{w}(\bv{s},t) := \bv{\zeta}(\bv{s},t) - \bv{s}$ ($\bv{s} \in B$), which is therefore related to the velocity as follows:
\begin{equation}
\label{eq: w u rel}
\dot{\bv{w}}(\bv{s},t) = \bv{u}(\bv{x},t)\big|_{\bv{x} = \bv{\zeta}(\bv{s},t)}.
\end{equation}
The practical enforcement of Eq.~\eqref{eq: w u rel} is crucial to distinguish an immersed method from another.  In the original \acro{IBM}, Eq.~\eqref{eq: w u rel} was enforced via Dirac-$\delta$ distributions.  In the present formulation, we enforce Eq.~\eqref{eq: w u rel} variationally. Our approach can be seen as a special case of what is discussed in~\cite{Boffi2014The-Finite-Elem-0}, where Eq.~\eqref{eq: w u rel} is enforced through a Lagrange multiplier.

\subsection{Functional setting}
\label{subsec: Functional setting}
The principal unknowns of our fluid-structure interaction problem are the fields
\begin{equation}
\label{eq: unknowns}
\bv{u}(\bv{x},t), \quad
p(\bv{x},t), \quad \text{and} \quad
\bv{w}(\bv{s},t),
\quad\text{with $\bv{x} \in \Omega$, $\bv{s} \in B$, and $t \in [0,T)$.}
\end{equation}
The functional spaces for these fields are
\begin{gather}
\label{eq: functional space u}
\bv{u} \in \mathscr{V} = H_{D}^{1}(\Omega)^{d} := \Bigl\{ \bv{u} \in L^{2}(\Omega)^{d} \,\big|\, \nabla_{\bv{x}} \bv{u} \in L^{2} (\Omega)^{d \times d},  \bv{u}|_{\partial\Omega_{D}} = \bv{u}_{g} \Bigr\},
\\
\label{eq: functional space p}
p \in \mathscr{Q} := L^{2}(\Omega), \\
\label{eq: functional space w}
\bv{w} \in \mathscr{Y} = H^{1}(B)^{d}  := \Bigl\{ \bv{w} \in L^{2}(B)^{d} \,\big|\, \nabla_{\bv{s}} \bv{w} \in L^{2} (B)^{d \times d}  \Bigr\},
\end{gather}
where $\nabla_{\bv{x}}$ and $\nabla_{\bv{s}}$ denote the gradient operators relative to $\bv{x}$ and $\bv{s}$, respectively.

For convenience, we will use a prime to denote partial differentiation with respect to time:
\begin{equation}
\label{eq: prime notation}
\bv{u}'(\bv{x},t) := \frac{\partial \bv{u}(\bv{x},t)}{\partial t}
\quad \text{and} \quad
\bv{w}'(\bv{s},t) := \frac{\partial \bv{w}(\bv{s},t)}{\partial t}.
\end{equation}
Equations~\eqref{eq: functional space u} and~\eqref{eq: functional space p} imply that the fields
$\bv{u}$ and $p$ are defined everywhere in $\Omega$.  Because $\bv{u}$
is defined everywhere in $\Omega$, the function $\bv{\sigma}^{v}_{\f}$ is defined everywhere in $\Omega$ as well. For consistency, we must also extend the domain of definition of the mass density of the fluid.  Hence, we formally assume that
\begin{equation}
\label{eq: rho f domain of definition}
\rho_{\f} \in L^{\infty}(\Omega).
\end{equation}
Referring to Eq.~\eqref{eq: functional space u}, the function space for the velocity test functions is  $\mathscr{V}_{0}$, defined as
\begin{equation}
\label{eq: space of test functions v}
\mathscr{V}_{0} = H_{0}^{1}(\Omega)^{d} := \Bigl\{ \bv{v} \in L^{2}(\Omega)^{d} \,\big|\, \nabla_{\bv{x}} \bv{v} \in L^{2} (\Omega)^{d \times d},  \bv{v}|_{\partial\Omega_{D}} = \bv{0} \Bigr\}.
\end{equation}

%

\subsection{Governing equations: incompressible solid}
\label{subsec: Governing equations: incompressible solid}
When the solid is incompressible, the mass densities the fluid and the solid are constant, and the governing equations can be given the following form \citep{Heltai2012Variational-Imp0}:
\begin{gather}
\label{eq: Bmomentum weak partitioned last}
\begin{multlined}[b]
\int_{\Omega} \rho_{\f} \bigl[\bv{u}' + (\nabla_{\bv{x}}\bv{u})\bv{u} - \bv{b}\bigr] \cdot \bv{v} \d{v}
- \int_{\Omega} p (\nabla_{\bv{x}} \cdot \bv{v}) \d{v}
+
\int_{\Omega} \bv{\sigma}^{v}_{\f} \cdot \nabla_{\bv{x}}\bv{v} \d{v}
-\int_{\partial\Omega_{N}} \mathfrak{s}_{g} \cdot \bv{v} \d{a}
\\
+ 
\int_{B} \bigl\{[\rho_{\s_{0}}(\bv{s}) - \rho_{\f} J(\bv{s},t)] \bigl\{\bv{u}'(\bv{x},t) + \bigl[\nabla_{\bv{x}}\bv{u}(\bv{x},t)\bigl] \bv{u}(\bv{x},t) - \bv{b}(\bv{x},t)\bigr\}\cdot \bv{v}(\bv{x})\bigr|_{\bv{x} = \map} \d{V}
\\
+
\int_{B} J(\bv{s},t) \bigl(\bv{\sigma}^{v}_{\s} - \bv{\sigma}^{v}_{\f}\bigr) \cdot \nabla_{\bv{x}}\bv{v}(\bv{x})\bigr|_{\bv{x} = \map} \d{V}
\\
+
\int_{B} \tensor{P}^{e}_{\s} \, \trans{\tensor{F}}(\bv{s},t) \cdot \nabla_{\bv{x}}\bv{v}(\bv{x})\bigr|_{\bv{x} = \map} \d{V}
=
0
\quad \forall \bv{v} \in \mathscr{V}_{0},
\end{multlined}
\\
\label{eq: Bmass weak partitioned}
\int_{\Omega} q (\nabla \cdot \bv{u}) \d{v} = 0
\quad \forall q \in \mathscr{Q},
\\
\label{eq: w u rel weak}
\Phi_{B}
\int_{B} \Bigl[\bv{w}'(\bv{s},t) - \bv{u}(\bv{x},t)\big|_{\bv{x} = \map}\Bigr] \cdot \bv{y}(\bv{s}) \d{V} = 0
\quad
\forall \bv{y} \in \mathscr{Y},
\end{gather}
where $\Phi_{B}$ is a constant with dimensions of mass over time divided by length cubed, i.e., dimensions such that, in 3D, the volume integral of the quantity $\Phi_{B}\bv{w}'$ has the same dimensions as a force.

Equation~\eqref{eq: Bmomentum weak partitioned last} is obtained from the second of Eqs.~\eqref{eq: Balance of mass and momentum} in a conventional way, i.e., by first constructing the scalar product with a test function $\bv{v} \in \mathscr{V}_{0}$, then integrating over $\Omega$, and finally applying the divergence theorem.  As the equation in question is supported over the entire domain $\Omega$, we proceed to treat integrals over $\Omega$ in the following manner.  Let $\phi$ represent a generic quantity with expressions $\phi_{\f}$ and $\phi_{\s}$ over the fluid and solid domains, respectively.  Also, let $\check\phi_{\f}$ represent the extension of $\phi_{\f}$ over $\Omega$ (see Eq.~\eqref{eq: rho f domain of definition} and discussion proceeding it). Then, we have
\begin{equation}
\label{eq: ints over Omega}
\begin{aligned}[b]
\int_{\Omega} \phi \d{v} &= \int_{\Omega\setminus B_{t}} \phi_{\f} + \int_{B_{t}} \phi_{\s} \d{v}
\\
&=\int_{\Omega} \check{\phi}_{\f} \d{v} + \int_{B_{t}} (\phi_{\s} - \check{\phi}_{\f}) \d{v}
\\
&= \int_{\Omega} \check{\phi}_{\f} \d{v} + \int_{B} [(\phi_{\s} - \check{\phi}_{f}) \circ \bv{\zeta}] J \d{V},
\end{aligned}
\end{equation}
where the last term in the above expression is simply the evaluation of the integral over $B_{t}$ as an integral over the reference configuration $B$.  

Equation~\eqref{eq: Bmass weak partitioned} is obtained from the first of Eq.~\eqref{eq: Balance of mass and momentum} by constructing the scalar product with test functions in $\mathscr{Q}$ (see Eq.~\eqref{eq: functional space p}) and then integrating over $\Omega$.  Finally, Eq.~\eqref{eq: w u rel weak} is obtained by constructing the scalar product of Eq.~\eqref{eq: w u rel} with test functions in $\mathscr{Y}$ (see Eq.~\eqref{eq: functional space w}) and integrating over $B$.

We note that a key element of any fully variational formulation of immersed methods is (the variational formulation of) the equation enabling the tracking of the motion of the solid, here Eq.~\eqref{eq: w u rel weak}. In the discrete formulation, this relation is as general as the choice of the finite-dimensional functional subspaces approximating $\mathscr{V}$ and $\mathscr{Y}$ and it is key to the stability of the method.  We note that a new approach for the solid motion equation has recently been presented by \cite{Boffi2014The-Finite-Elem-0} employing distributed Lagrange multipliers. As it turns out, this formulation, which can handle problems with variable density and viscosity, yields unconditionally stable semi-implicit time advancing scheme.

\subsection{Governing equations: compressible solid}
\label{subsec: Governing equations: compressible solid}
When the solid is compressible, the contribution of $p$ to the balance of linear momentum and the incompressibility constraint must be restricted to the domain $\Omega\setminus B_{t}$. Therefore, the weak problem is rewritten as follows:
\begin{gather}
\label{eq: Bmomentum weak partitioned last compressible solid}
\begin{multlined}[b]
\int_{\Omega} \rho_{\f} (\dot{\bv{u}} - \bv{b}) \cdot \bv{v} \d{v}
- \int_{\Omega} p \nabla\cdot \bv{v} \d{v}
+
\int_{\Omega} \bv{\sigma}^{v}_{\f} \cdot \nabla_{\bv{x}}\bv{v} \d{v}
-\int_{\partial\Omega_{N}} \mathfrak{s}_{g} \cdot \bv{v} \d{a}
\\
+ 
\int_{B} \bigl\{[\rho_{\s_{0}}(\bv{s}) - \rho_{\f_0}] [\dot{\bv{u}}(\bv{x},t) - \bv{b}(\bv{x},t)]\cdot \bv{v}(\bv{x})\bigr|_{\bv{x} = \map} \d{V}
\\
+ \int_{B} J(\bv{s},t) p(\bv{x},t) \nabla\cdot \bv{v}(\bv{x})\bigr|_{\bv{x} = \map} \d{V}
\\
+
\int_{B} J(\bv{s},t) \bigl(\bv{\sigma}^{v}_{\s} - \bv{\sigma}^{v}_{\f}\bigr) \cdot \nabla_{\bv{x}}\bv{v}(\bv{x})\bigr|_{\bv{x} = \map} \d{V}
\\
+
\int_{B} \tensor{P}^{e}_{\s} \, \trans{\tensor{F}}(\bv{s},t) \cdot \nabla_{\bv{x}}\bv{v}(\bv{x})\bigr|_{\bv{x} = \map} \d{V}
=
0
\quad \forall \bv{v} \in \mathscr{V}_{0}.
\end{multlined}
\\
\label{eq: balance of mass restricted}
\int_{\Omega} q \nabla\cdot \bv{u} \d{v} - \int_{B_{t}} q \nabla\cdot \bv{u} \d{v} = 0,
\\
\label{eq: w u rel weak compressible solid}
\Phi_{B} \int_{B} \Bigl[\dot{\bv{w}}(\bv{s},t) - \bv{u}(\bv{x},t)\big|_{\bv{x} = \map}\Bigr] \cdot \bv{y}(\bv{s}) \d{V} = 0
\quad
\forall \bv{y} \in \mathscr{Y}.
\end{gather}
Equations~\eqref{eq: Bmomentum weak partitioned last compressible solid}--\eqref{eq: w u rel weak compressible solid} would allow us to determine a unique solution if the field $p$ were restricted to the domain $\Omega\setminus B_{t}$.  However, our numerical scheme  requires the domain of $p$ to be $\Omega$.  To sufficiently constraint the behavior of $p$ over $B_{t}$, we have considered two possible strategies.  One is to set to zero the restriction of $p$ to $B_{t}$.  The other is physically motivated and based on the fact that, for a Newtonian fluid, $p$ is the mean normal stress.  Hence, $p$ can be constrained on $B_{t}$ so to represent the mean normal stress everywhere in $\Omega$. Given the constitutive response function of a compressible solid material, the corresponding mean normal stress is
\begin{equation}
\label{eq: Mean normal stress def}
p_{\s}[\bv{u},\bv{w}] = -\frac{1}{\trace{\tensor{I}}}
\Bigr[ \bv{\sigma}_{\s}^{v}[\bv{u}] \cdot \tensor{I} + J^{-1}[\bv{w}] \tensor{P}_{\s}^{e}[\bv{w}] \cdot \tensor{F}[\bv{w}]\Bigl].
\end{equation}
With this in mind, we replace Eq.~\eqref{eq: balance of mass restricted} with the following equation:
\begin{multline}
\label{eq: Bmass weak partitioned compressible solid}
-\int_{\Omega} q \nabla\cdot \bv{u} \d{v} 
+ \int_{B} J(\bv{s},t) q(\bv{x}) \nabla\cdot \bv{u}(\bv{x},t)\bigr|_{\bv{x} = \map} \d{V} 
\\
+ \int_{B} c_{1}
J(\bv{s},t) \bigr[ p(\bv{x},t) - c_{2} p_{\s}[\bv{u},\bv{w}] \bigl]
q(\bv{x})\bigr|_{\bv{x} = \map} \d{V} = 0
\quad \forall q \in \mathscr{Q},
\end{multline}
where $c_{1} > 0$ is a constant parameter with dimensions of length times mass over time, and where $c_{2}$ is a dimensionless constant with values $0$ or $1$.  For $c_{2} = 0$, $p = 0$ (weakly) over $B_{t}$, whereas for $c_{2} = 1$, $p$ is (weakly) constrained to be the mean normal stress everywhere in $\Omega$.

Equations presented above can be compactly reformulated in terms of the Hilbert space $\mathscr{Z} := \mathscr{V}\times \mathscr{Q}\times \mathscr{Y}$, and $\mathscr{Z}_{0} := \mathscr{V}_{0} \times \mathscr{Q} \times \mathscr{Y}$ with inner product given by the sum of the inner products of the generating spaces. Defining $\mathscr{Z} \ni \xi := \trans{[\bv{u}, p, \bv{w}]}$ and $\mathscr{Z}_{0} \ni \psi := \trans{[\bv{v}, q, \bv{y}]}$, we can given our problem the following form:
\begin{problem}[Grouped dual formulation]
\label{prob: IFG}
Given an initial condition $\xi_0 \in \mathscr{Z}$, for all $t \in (0,T)$ find $\xi(t) \in \mathscr{Z}$, such that
\begin{equation}
\label{eq:formal grouped dual}
\langle \mathcal{F}(t, \xi, \xi') , \psi \rangle =0, \quad \forall \psi \in \mathscr{Z}_0,
\end{equation}
where $\mathcal{F} : \mathscr{Z} \mapsto \mathscr{Z}_0^*$, $\mathscr{Z}_0^*$ is the dual of $\mathscr{Z}_0$, and the precise definition of the operator $\mathcal{F}$ can be deduced from the integral expressions in Eqs.~\eqref{eq: Bmomentum weak partitioned last}--\eqref{eq: w u rel weak} for the incompressible case and Eqs.~\eqref{eq: Bmomentum weak partitioned last compressible solid}, \eqref{eq: w u rel weak compressible solid}, and~\eqref{eq: Bmass weak partitioned compressible solid} for the compressible case (see \citealp{Heltai2012Variational-Imp0} for additional details).
\end{problem}
\begin{remark}[Initial condition for the pressure]
  In Problem~\ref{prob: IFG}, an initial condition for the triple $\xi_0 = \trans{[\bv{u}_0, p_0, \bv{w}_0]}$ is required just as a matter of compact representation of the problem. However, only the initial conditions  $\bv{u}_0$ and $\bv{w}_0$ are used, since there is no time derivative acting on $p$.
\end{remark}

\begin{remark}[Energy estimates]
\cite{Heltai2012Variational-Imp0} have shown that the formulation presented thus far leads to energy estimates that that are formally identical to those of the continuous formulation.
\end{remark}

\subsection{Spatial Discretization by finite elements}
\label{section: discretization by FEM}
Domains $\Omega$ and $B$ are decomposed into independent triangulations $\Omega_{h}$ and $B_{h}$, respectively, consisting of cells $K$ (triangles or quadrilaterals in 2D, and tetrahedra or hexahedra in 3D) such that
\begin{enumerate}
\item
$\overline{\Omega} = \cup \{ K \in \Omega_{h} \}$, and $\overline{B} = \cup \{ K \in B_{h} \}$;

\item
Any two cells $K,K'$ only intersect in common faces, edges, or vertices;

\item
The decomposition $\Omega_{h}$ matches the decomposition $\partial \Omega = \partial\Omega_{D} \cup \partial\Omega_{N}$.
\end{enumerate}
On $\Omega_{h}$ and $B_{h}$, finite dimensional subspaces $\mathscr{V}_{h} \subset \mathscr{V}$, $\mathscr{Q}_{h} \subset \mathscr{Q}$, and $\mathscr{Y}_{h} \subset \mathscr{Y}$ are defined such that
\begin{alignat}{5}
\label{eq: functional space u h}
\mathscr{V}_h &:= \Bigl\{ \bv{u}_h \in \mathscr{V} \,&&\big|\, \bv{u}_{h|K}  &&\in \mathcal{P}_V(K), \, K &&\in \Omega_h \Bigr\} &&\equiv \vssp\{ \bv{v}_{h}^{i} \}_{i=1}^{N_{V}},
\\
\label{eq: functional space p h}
\mathscr{Q}_h &:= \Bigl\{ p_h \in \mathscr{Q} \,&&\big|\, p_{h|K}  &&\in \mathcal{P}_Q(K), \, K &&\in \Omega_h \Bigr\} &&\equiv \vssp\{ q_{h}^{i} \}_{i=1}^{N_{Q}},
\\
\label{eq: functional space w h}
\mathscr{Y}_{h} &:= \Bigl\{ \bv{w}_h \in \mathscr{Y} \,&&\big|\, \bv{w}_{h|K} &&\in \mathcal{P}_Y(K), \, K &&\in B_h \Bigr\} &&\equiv \vssp\{ \bv{y}_{h}^{i} \}_{i=1}^{N_{Y}},
\end{alignat}
where $\mathcal{P}_{V}(K)$, $\mathcal{P}_{Q}(K)$ and $\mathcal{P}_{Y}(K)$ are polynomial spaces of degree $r_{V}$, $r_{Q}$ and $r_{Y}$ respectively on the cells $K$, and $N_V$, $N_Q$ and $N_Y$ are the dimensions of each finite dimensional space.  Notice that we use only one discrete space for both $\bv{u}$ and $\bv{u}'$ and one for $\bv{w}$ and $\bv{w}'$, even though the continuous functional spaces should be different, as discussed in details in~\cite{Heltai2012Variational-Imp0}.  Also, we chose the pair $\mathscr{V}_{h}$ and $\mathscr{Q}_{h}$ so as to satisfy the inf-sup condition for existence, uniqueness, and stability of solutions of the Navier-Stokes component of the problem (see, e.g., \citealp{BrezziFortin-1991-a}).

\subsection{Variational velocity coupling}
\label{subsection: FEM implementation}
The implementation of operator with support over $\Omega$ is common in \acro{FEM} approaches to the Navier-Stokes equation.  What is less common in the practical implementation of operators defined over $B_{h}$ involving the evaluation of fields supported over $\Omega_{h}$.  Here we outline the basic features of the implementation (see also \citealp{Heltai2012Variational-Imp0}).  As a specific example, we describe the treatment of the operator coupling the velocities of the fluid and of the solid domain.  Referring to the terms involving $\Phi_{B} \bv{u}(\bv{x},t)\big|_{\bv{x} = \map}\cdot \bv{y}(\bv{s})$ in Eqs.~\eqref{eq: w u rel weak} and~\eqref{eq: w u rel weak compressible solid}, we define a matrix $M_{wu}(\bv{w})$ whose $ij$ element is given by
\begin{equation}
  \label{eq: MGamma def bis}
  {M}_{wu}^{ij}(\bv{w}_{h}) 
  = \Phi_{B} \int_{B}  \bv{v}^j_{h}(\bv{x})\big|_{\bv{x} = \bv{s} + \bv{w}_{h}(\bv{s},t)} \cdot \bv{y}^i_{h}(\bv{s}) \d{V}.
\end{equation}
The computation of the integral in Eq.~\eqref{eq: MGamma def bis} is done by summing the contributions due to each cell $K$ in $B_{h}$ via a quadrature rules with $N_Q$ points. The functions $\bv{y}^i_{h}(\bv{s})$ have support over $B_{h}$, whereas the functions $\bv{v}^j_{h}(\bv{x})$ (with $\bv{x} = \bv{s} + \bv{w}_{h}(\bv{s},t)$) have support on $\Omega_{h}$.
\begin{figure}[htb]
    \centering
    \includegraphics{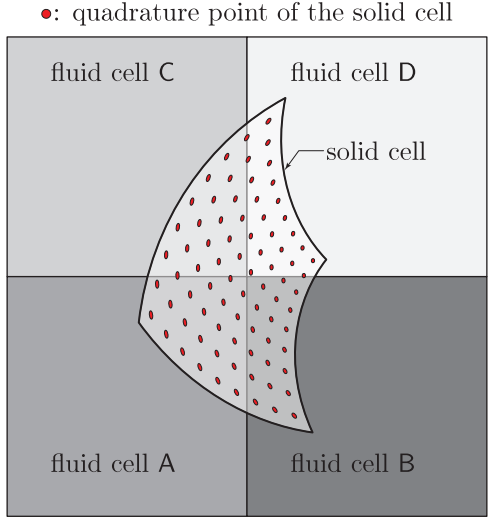}
    \caption{Cells denote as \textsf{A}--\textsf{D} represent a four-cell patch of the triangulation of the fluid domain.  The cell denoted as ``solid cell'' represents a cell of the triangulation of the immersed solid domain that is contained in the union of cells \textsf{A}--\textsf{D} of the fluid domain.  The filled dots represent the quadrature points of the quadrature rule adopted to carry out integration over the cells of the immersed domain.}
    \label{fig: integration}
\end{figure}
Hence, we first determine the position of the quadrature points of the solid element, both relative to the reference unit element and relative to the underlying global coordinate system through the mappings:
\begin{alignat}{3}
  \label{eq:mapping Khat K solid}
  \bv{s}_K & : \hat{K} := [0,1]^d &&\mapsto K \in B_h, \\
  \label{eq:mapping K K solid}
  I+\bv{w}_h & : K && \mapsto K_t \in B_{t,h}.
\end{alignat}
Then, these global coordinates are passed to an algorithm that finds the fluid cells in $\Omega_{h}$ containing the points in question.  The outcome of this operation is sketched in Fig.~\ref{fig: integration} where we show the deformed configuration ($K_t$ in Eq.~(\ref{eq:mapping K K solid})) of a cell of $B_{h}$ ($K$ in Eq.~(\ref{eq:mapping Khat K solid})) straddling four cells of $\Omega_{h}$ denoted fluid cells \textsf{A}--\textsf{D}.  The quadrature points over the solid cell are represented by the filled circles.

\subsection{Time discretization}
\label{sec:time-discretization}
The discrete counterpart of Eq.~\eqref{eq:formal grouped dual} represents a system of nonlinear differential algebraic equations (\acro{DAE}).  The time derivative $\xi'$ is approximated very simply via an implicit-Euler scheme:
  \begin{equation}
    \label{eq:time derivative theta}
    \xi_n' = h^{-1} \bigl( \xi_{n} - \xi_{n-1} \bigr),
  \end{equation}
where $\xi_{n}$ and $\xi_{n}'$ are the computed approximations to $\xi(t_{n})$ and $\xi'(t_{n})$, respectively, and the step size $h = t_{n} - t_{n-1}$ is kept constant throughout the computation.  Although not second order accurate, this time stepping scheme is asymptotically stable.  Overall, our system of discretized equations is solved in a fully coupled implicit manner.  A more comprehensive discussion of other fully implicit methods has been presented in \cite{Heltai2012A-Fully-Coupled0}.  Other methods, not necessarily implicit are possible but have not been explored yet, at least by the authors. As noted earlier, \cite{Boffi2014The-Finite-Elem-0} have recently presented a variational formulation using distributed Lagrange multipliers, with an interesting unconditionally stable time advancing scheme.

Going back to the discussion of our proposed methodology, the advancement from a time $t^{k}$ to $t^{k+1}$ requires the solution of a Newton iteration cycle, which, in progressing from iterate $n$ to iterate $n+1$, takes on the following matrix form:
\begin{equation}
\label{eq: matrix form newton step}
\begin{pmatrix}
\mathsf{A}^{n+1} & \mathsf{B}^{\mathrm{T}} & \mathsf{A}_{\s}^{n+1}
\\
\mathsf{B} & 0 & 0
\\
-M_{wu}^{n+1} & 0 & M_{ww}
\end{pmatrix}
\begin{pmatrix}
\delta\bv{u}_{h}
\\
\delta p_{h}
\\
\delta \bv{w}_{h}
\end{pmatrix}
=
\begin{pmatrix}
\bv{R}_{u}^{n+1}
\\
R_{p}^{n+1}
\\
\bv{R}_{w}^{n+1}
\end{pmatrix},
\end{equation}
where, cognizant that the elements in the above equation depend on the input data at time $t^{k+1}$ and the solution at time $t^{k}$,
\begin{enumerate}[(i)]
\item
$\delta \bv{u}_{h} = \bv{u}_{h}^{n+1} - \bv{u}_{h}^{n}$, $\delta p_{h} = p_{h}^{n+1} - p_{h}^{n}$, $\delta \bv{w}_{h} = \bv{w}_{h}^{n+1} - \bv{w}_{h}^{n}$;

\item
$\bv{R}_{u}$, $R_{p}$, and $R_{w}$ are the residuals for equations Eq.~\eqref{eq: Bmomentum weak partitioned last}, \eqref{eq: Bmass weak partitioned}, and \eqref{eq: w u rel weak}, respectively;
\item
blocks $(1,1)$, $(1,2)$, $(2,1)$, and $(2,2)$ are precisely those of a pure Navier-Stokes problem;

\item
$\mathsf{A}_{\s}$ describes set of forces on the fluid due to the solid's response;

\item
$M_{ww}$ is the mass matrix associated to the field $\bv{w}$ over the solid's domain, and $M_{wu}$ is the operator discussed in Eq.~\eqref{eq: MGamma def bis}.
\end{enumerate}
In the compressible case one obtains a similar problem except for the fact that the $(2,2)$ block is not identically zero due to the fact that the pressure behavior over the domain $B_{t}$ is being constrained as indicated by the last term in Eq.~\eqref{eq: Bmass weak partitioned compressible solid}.

\section{Numerics}
\label{sec: Numerics}
To the authors' knowledge, immersed methods have not been validated as many \acro{ALE} methods have been via the rigorous benchmark tests by \citealp{Turek2006Proposal-for-Nu0} due to intrinsic modeling restrictions.  In this paper we present results that show that the method in \cite{Heltai2012Variational-Imp0}, being applicable to the physical systems in the benchmarks tests by \citealp{Turek2006Proposal-for-Nu0}, can indeed satisfy these benchmarks in a rigorous way.  As such, to the best of the authors' knowledge, the results shown herein are the first in which an immersed method is shown to satisfy the benchmarks in question.  We point out that in previous works concerning the use of fully variational approaches to the immersed finite element method, the immersed body was assumed to be incompressible and viscoelastic (with a linear viscous component formally similar to that of the fluid) (cf.\ \citealp{Heltai_2006_The-Finite_0,Heltai-2008-a,BoffiGastaldiHeltaiPeskin-2008-a,Heltai2012Variational-Imp0}).  In this paper, we present for the first time results in which the immersed body is compressible and purely hyperelastic. Furthermore, we consider cases in which the immersed solid and the fluid have different densities as well as dynamic viscosities (when the solid is assumed to have a linear viscous component to its stress response). We want to emphasize that in all simulations, even those with a compressible elastic body, the fluid is always modeled as incompressible, as opposed to nearly incompressible.

All the results presented in this section have been obtained using a modified version of the code presented in~\cite{Heltai2012A-Fully-Coupled0}, implemented using the \texttt{deal.II} library (see, e.g.,~\citealp{BangerthHartmannKanschat-2007-a,BangerthHeisterHeltai-2015-a}).

\subsection{Discretization}
The approximation spaces we used in our simulations for the approximations of the velocity field $\bv{u}_h$ and of the displacement field $\bv{w}_h$ are the piecewise bi-quadratic spaces of continuous vector functions over $\Omega$ and over $B$, respectively, which we will denote by $\mathcal{Q}^2_{0}$ space.\footnote{In general, we denote by $\mathcal{Q}_{c}^{p}$ the space of piecewise polynomials consisting of tensor products of polynomials of order $p$ and with global continuity degree $c$.  We denote by $\mathcal{P}_{c}^{p}$ the spaces of piecewise polynomials of maximum order $p$ and continuity degree $c$.  In all cases a subscript $c = -1$ denotes discontinuous spaces.} For the pressure field $p$, in some cases we have used the piecewise continuous bi-linear space $\mathcal{Q}_{0}^{1}$ and in other cases the piecewise discontinuous linear space $\mathcal{P}^1_{-1}$ over $\Omega$.  Both the $\mathcal{Q}_{0}^2\text{-}\vert\mathcal{Q}_{0}^1$ and the $\mathcal{Q}_{0}^2\vert\mathcal{P}_{-1}^1$ pairs of spaces are known to satisfy the inf-sup condition for the approximation of the Navier-Stokes part of our equations (see, e.g., \citealp{BrezziFortin-1991-a}). The choice of the space $\mathcal{Q}_{0}^2$ for the displacement variable $\bv{w}_h$ is a natural choice, given the underlying velocity field $\bv{u}_h$. With this choice of spaces, Eqs.~\eqref{eq: w u rel weak} and~\eqref{eq: w u rel weak compressible solid} can be satisfied exactly when the solid and the fluid meshes are matching. 

\subsection{Results for Incompressible Immersed Solids}

\subsubsection{Static equilibrium of an annular solid comprising circumferential fibers and immersed in a stationary fluid}
This numerical test is motivated by the ones presented in \cite{BoffiGastaldiHeltaiPeskin-2008-a,GriffithLuo-2012-a} and it pertains to the case of an incompressible solid with both elastic and viscous components of the stress response. The objective of this test is to compute the equilibrium state of an initially undeformed thick annular cylinder submerged in a stationary incompressible fluid that is contained in a rigid prismatic box having a square cross-section.

The simulation is two-dimensional and comprises an annular solid with inner radius $R$ and thickness $w$, and filled with a stationary fluid that is contained in a square box of edge length $l$ (see Fig.~\ref{fig:RingEqm-Geometry}).  The reference and  deformed configurations can be described via polar coordinate systems with origins at the center of the annulus and whose unit vectors are given by $\left(\hat{\bv{u}}_{R}, \hat{\bv{u}}_{\Theta} \right)$ and $\left(\hat{\bv{u}}_{r}, \hat{\bv{u}}_{\theta} \right)$, respectively.
\begin{figure}[htbp]
	\begin{center}
		\includegraphics{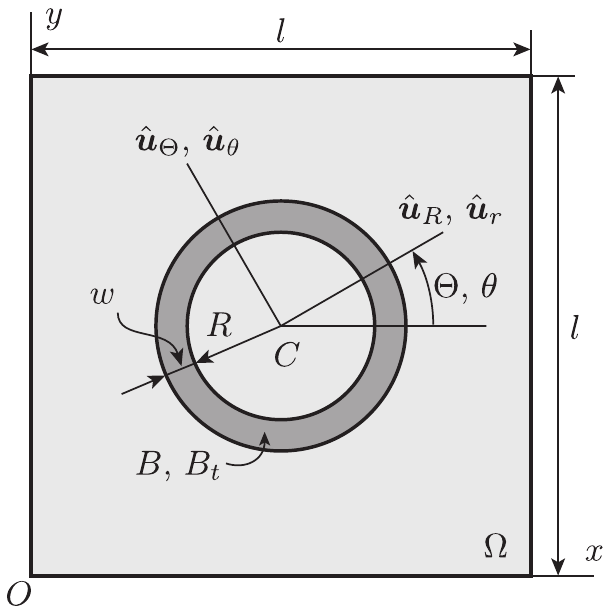}
	\caption{The reference and deformed configurations of a ring immersed in a square box filled with stationary fluid.}
	\label{fig:RingEqm-Geometry}
	\end{center}
\end{figure}
This ring is subjected to the hydrostatic pressure of the fluid $p_{i}$ and $p_{o}$ at its inner and outer walls, respectively. Negligible body forces act on the system and there is no inflow or outflow of fluid across the walls of the box. Since both the solid and the fluid are incompressible, neither the annulus nor the fluid will move and the problem reduces to determining the Lagrange multiplier field $p$. The elastic behavior of the ring is governed by a continuous distribution of concentric fibers lying in the circumferential direction. The first Piola-Kirchhoff and the Cauchy stress tensors are then given by, respectively,
\begin{equation}
\label{eq: PK stress for ring}
\tensor{P} = -p_{\s} \tensor{F}^{-T} + G^{e} \tensor{F} \hat{\bv{u}}_{\Theta} \otimes \hat{\bv{u}}_{\Theta}
\quad \text{and} \quad
\bv{\sigma}_{\s} = -p_{\s} \tensor{I} + G^{e} \hat{\bv{u}}_{\theta} \otimes \hat{\bv{u}}_{\theta},
\end{equation}
where $G^{e}$ is a constant modulus of elasticity, $p_{s}$ is the Lagrange multiplier that enforces incompressibility of the ring and $\hat{\bv{u}}_{\theta} = \hat{\bv{u}}_{\Theta}$ (since deformed and reference configurations coincide).  Recall that in the proposed immersed \acro{FEM} we have a single field $p$ representing the Lagrange multiplier everywhere, whether in the fluid or in the solid.  Therefore, we have $p = p_{\s}$ in the solid.  With this in mind, we observe that the equilibrium stress state in the fluid is purely hydrostatic.  Furthermore, since the boundary conditions on $\partial\Omega$ is of homogeneous Dirichlet type, the solution for the Lagrange multiplier $p$ over $\Omega$ is not unique.  We remove the non uniqueness by enforcing a zero average constraint on the field $p$.  Then it can be shown (cf.\ \citealp{Heltai2012A-Fully-Coupled0}) that the solution for the field $p$ is as follows:
\begin{equation}
p=
\begin{cases}
p_{o}=-\frac{\pi G^{e}}{2 l^{2}} \left( \left(R+w\right)^{2}-R^{2}\right)  & \mathrm{for} \quad R + w \leq r,\\
p_{s}= G^{e} \ln (\frac{R+w}{r})-\frac{\pi G^{e}}{2 l^{2}} \left( \left(R+w\right)^{2}-R^{2}\right) & \mathrm{for} \quad R <r < R+w,\\
p_{i}=G^{e} \ln (1+\frac{w}{R})-\frac{\pi G^{e}}{2 l^{2}} \left( \left(R+w\right)^{2}-R^{2}\right) &\mathrm{for} \quad  r \leq R,
\end{cases}
\label{eqn: p for ring eqm}
\end{equation}
with velocity of fluid $\bv{u}=\bv{0}$ and the displacement of the solid $\bv{w}=\bv{0}$. Note that Eq.~\eqref{eqn: p for ring eqm} is different from Eq.~(69) of \cite{BoffiGastaldiHeltaiPeskin-2008-a}, where $p$ varies linearly with $r$ (we believe this to be in error).

For all our numerical simulations we have used $R =\np[m]{0.25}$, $w=\np[m]{0.06250}$, $l=\np[m]{1.0}$ and $G^{e}=\np[Pa]{1}$ and for these values we obtain $p_{i}=\np[Pa]{0.16792}$ and $p_{o}=\np[Pa]{-0.05522}$ using Eq.~\eqref{eqn: p for ring eqm}. We have used $\rho=\np[kg/m^{3}]{1.0}$, dynamic viscosities $\mu_{\f} = \mu_{\s} = \mu=\np[Pa \!\cdot\! s]{1.0}$, and time step size $h=\np[s]{1e-3}$ in our tests. For all our numerical tests we have used $\mathcal{Q}_{0}^{2}$ elements to represent $\bv{w}$ of the solid, whereas we have used (i) $\mathcal{Q}_{0}^{2}\vert\mathcal{P}_{-1}^{1}$ elements, and (ii) $\mathcal{Q}_{0}^{2}\vert\mathcal{Q}_{0}^{1}$ elements to represent $\bv{v}$ and $p$ over the control volume. We present a sample profile of $p$ over the entire control volume and its variation along different values of $y$, after one time step, in Fig.~\ref{fig:DGP_Pressure} and Fig.~\ref{fig:FEQ_Pressure} for $\mathcal{Q}_{0}^{2}\vert\mathcal{P}_{-1}^{1}$ and $\mathcal{Q}_{0}^{2}\vert\mathcal{Q}_{0}^{1}$ elements, respectively.

The convergence rate  (see, Tables~\ref{tab:ring-eqm-dgp} and \ref{tab:ring-eqm-feq} for $\mathcal{Q}_{0}^{2}\vert\mathcal{P}_{-1}^{1}$ and $\mathcal{Q}_{0}^{2}\vert\mathcal{Q}_{0}^{1}$ elements, respectively) is 2.5 for the $L^{2}$ norm of the velocity, 1.5 for the $H^{1}$ norm of the velocity and 1.5 for the $L^{2}$ norm of the pressure which matches the rates presented in \cite{BoffiGastaldiHeltaiPeskin-2008-a}. In all these numerical tests we have used \np{1856}~cells with \np{15776}~DoFs for the solid.


\begin{figure}[htbp]
	\begin{center}
	\subfigure[Over the entire domain]
	{\includegraphics[width=3in]{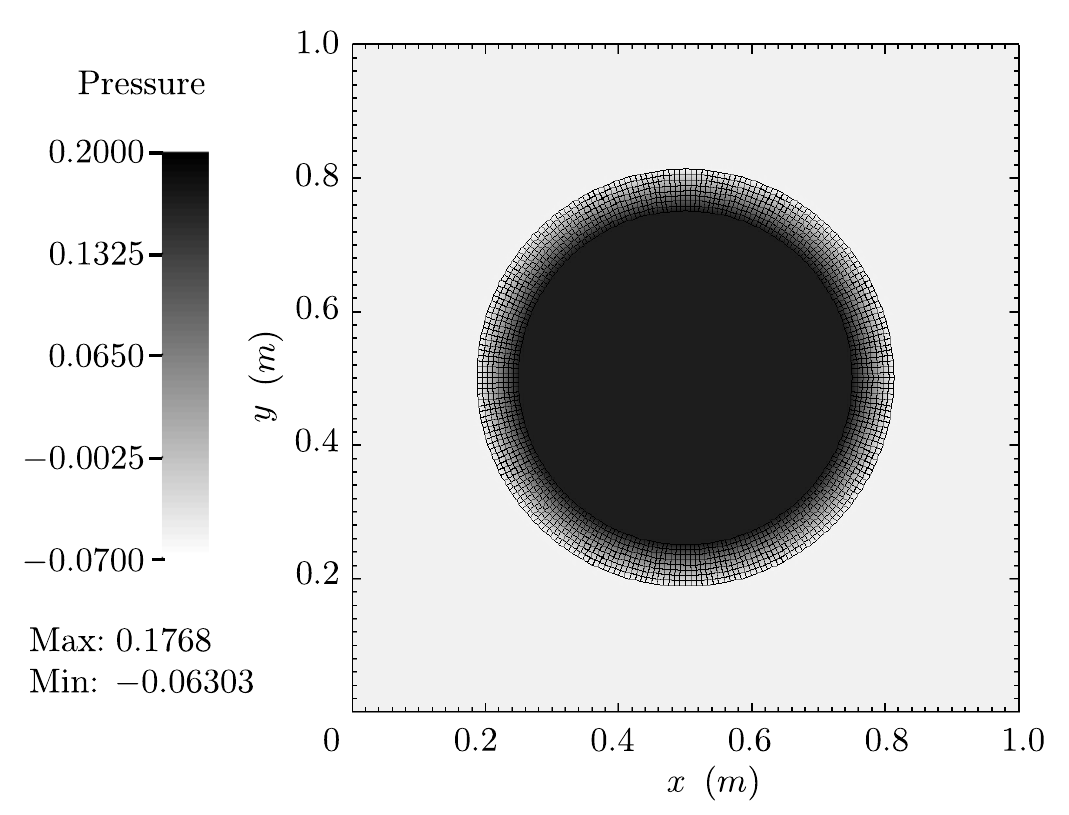}
	}
	\subfigure[At different values of $y$]
	{\includegraphics[width=3.1in]{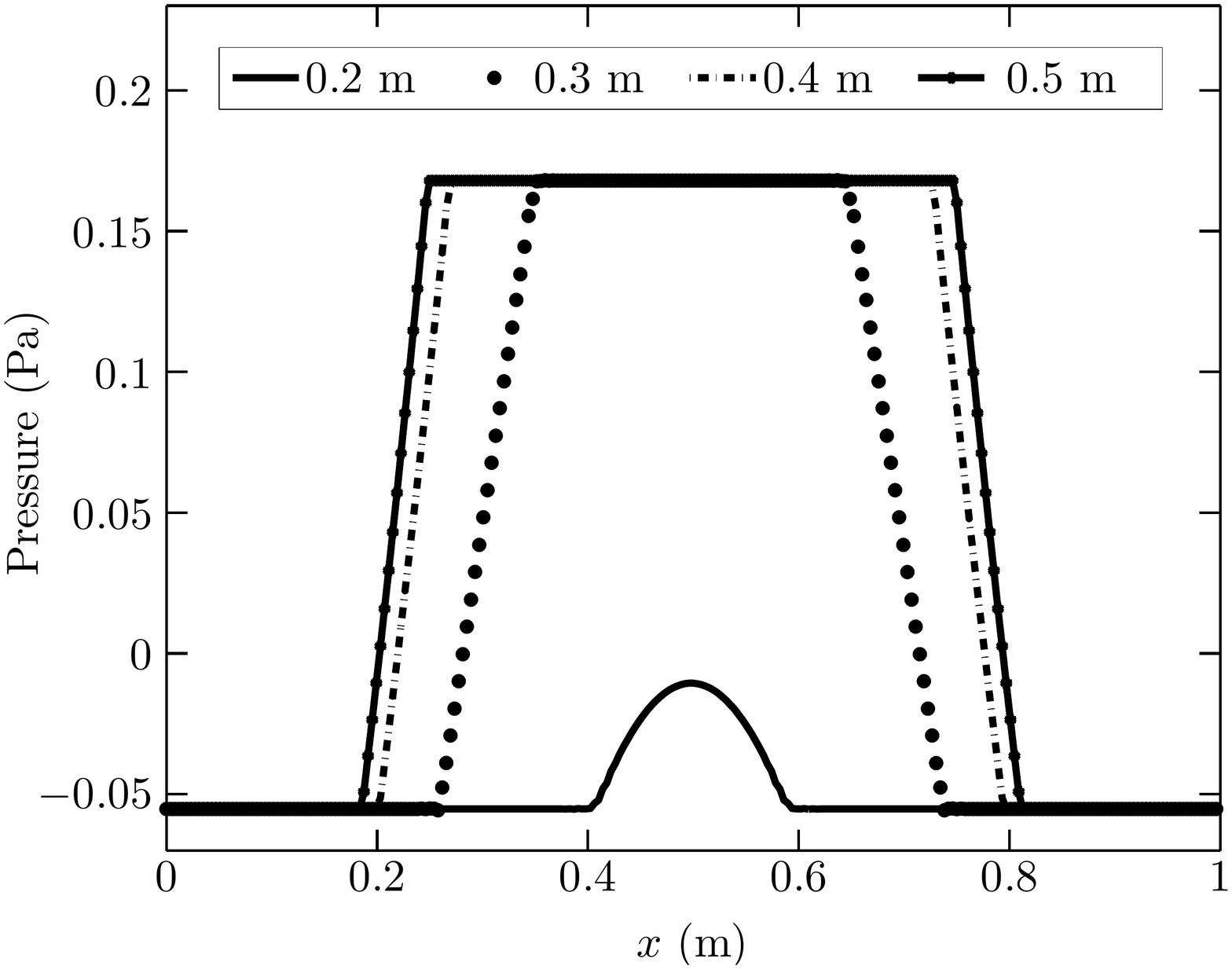}
	}
	\caption{The values of $p$ after one time step when using $\mathcal{P}_{-1}^{1}$ elements for $p$.}
	\label{fig:DGP_Pressure}
	\end{center}
\end{figure}
\begin{table}[htbp]\small
\caption{Error convergence rate obtained when using $\mathcal{P}_{-1}^{1}$ element for $p$ after one time step.}
\begin{center}
\begin{tabular*} {\textwidth} {@{\extracolsep{\fill}} c c c c c c c c}
\toprule
No. of cells &
No. of DoFs  &
$\|\mathbf{u}_{h}-\mathbf{u}\|_{0}$ & &
$\|\mathbf{u}_{h}-\mathbf{u}\|_{1}$ & &
$\|p_{h}-p\|_{0}$ & 
\\
\midrule
 \np{256} &   \np{2946}  & 2.00605e-05 & -    & 1.95854e-03 & -    & 6.71603e-03 & -    \\
\np{1024} &  \np{11522}  & 3.69389e-06 & 2.44 & 7.44696e-04 & 1.40 & 2.47476e-03 & 1.44 \\
\np{4096} &  \np{45570}  & 5.76710e-07 & 2.68 & 2.25134e-04 & 1.73 & 8.74728e-04 & 1.50 \\
\np{16384} & \np{181250} & 1.06127e-07 & 2.44 & 8.24609e-05 & 1.45 & 3.14028e-04 & 1.48 \\ 
\bottomrule
\end{tabular*}
\end{center}
\label{tab:ring-eqm-dgp}
\end{table}
\begin{figure}[htbp]
	\begin{center}
	\subfigure[Over the entire domain]
	{\includegraphics[width=3in]{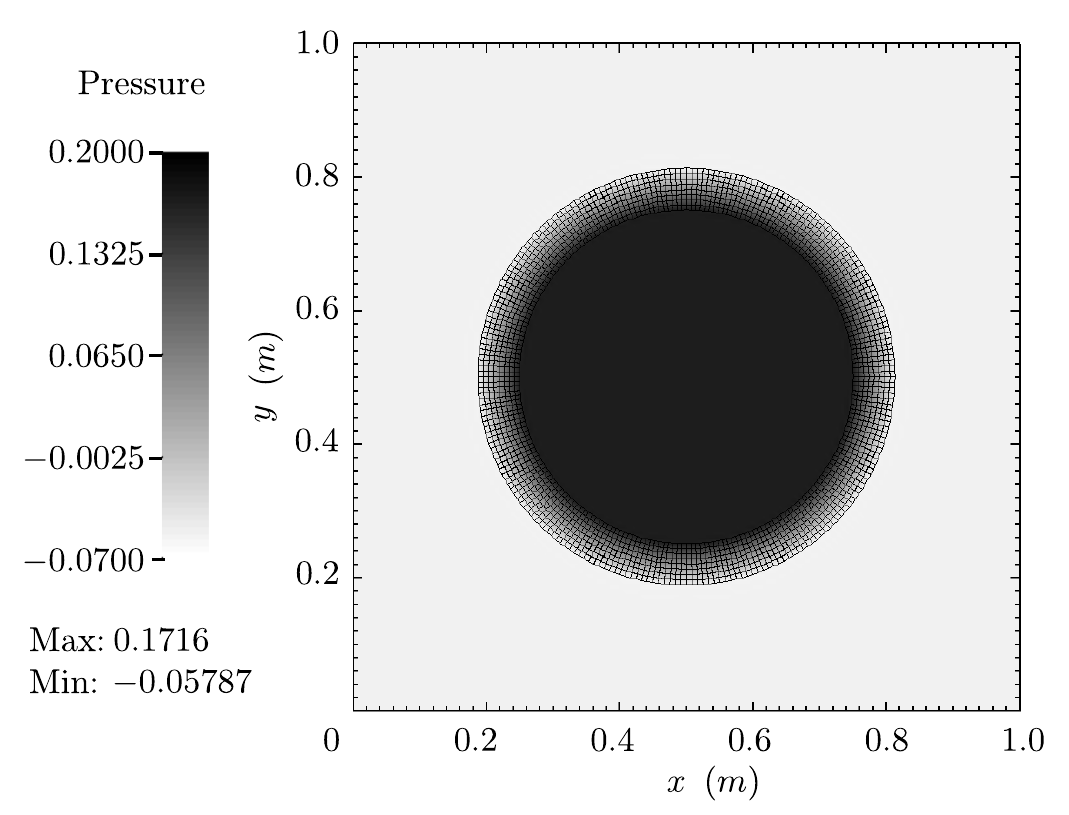}
	}
	\subfigure[At different values of $y$]
	{\includegraphics[width=3.1in]{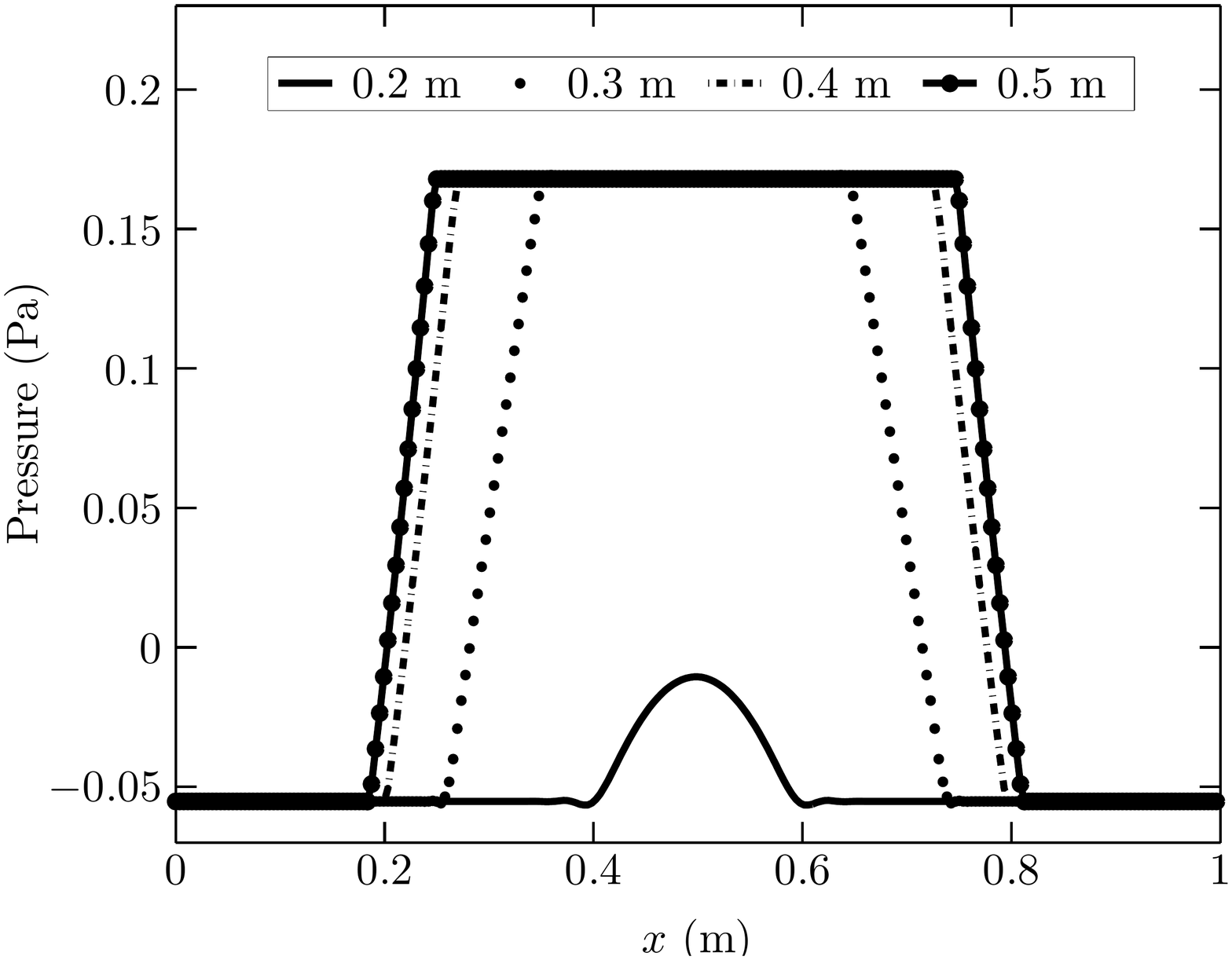}
	}
	\caption{The values of $p$ after one time step when using $\mathcal{Q}_{0}^{1}$ elements for $p$.}
	\label{fig:FEQ_Pressure}
	\end{center}
\end{figure}
\begin{table}[htbp]\small
\caption{Error convergence rate obtained when using $\mathcal{Q}_{0}^{1}$ element for $p$ after one time step.}
\begin{center}
\begin{tabular*} {0.95\textwidth} {@{\extracolsep{\fill}} c c c c c c c c}
\toprule
No. of cells &
No. of DoFs  &
$\|\mathbf{u}_{h}-\mathbf{u}\|_{0}$ & &
$\|\mathbf{u}_{h}-\mathbf{u}\|_{1}$ & &
$\|p_{h}-p\|_{0}$ & 
\\
\midrule
 \np{256} &   \np{2467} & 4.36912e-05 &    - & 2.79237e-03 &    - & 7.39310e-03 &- \\ 
\np{1024} &   \np{9539} & 6.14959e-06 & 2.83 & 9.02397e-04 & 1.63 & 2.42394e-03 & 1.61 \\ 
\np{4096} &  \np{37507} & 1.28224e-06 & 2.26 & 3.49329e-04 & 1.37 & 9.10608e-04 & 1.41 \\ 
\np{16384} & \np{148739} & 2.33819e-07 & 2.46& 1.25626e-04 & 1.48 & 3.27256e-04 & 1.48 \\
\bottomrule
\end{tabular*}
\end{center}
\label{tab:ring-eqm-feq}
\end{table}

\subsubsection{Disk entrained in a lid-driven cavity flow}
We test the volume conservation of our numerical method by measuring the change in the area of a disk that is entrained in a lid-driven cavity flow of an incompressible, linearly viscous fluid. This is another example pertaining to a system with an incompressible solid with stress response that includes both elastic and viscous contributions. This test is motivated by similar ones presented in \cite{WangZhang_2010_Interpolation-functions-0,GriffithLuo-2012-a}. Referring to Fig.~\ref{fig:LDCFlow-Ball-Geometry}, the disk has a radius $R=\np[m]{0.2}$ and its center $C$ is initially positioned at $x=\np[m]{0.6}$ and $y=\np[m]{0.5}$ in the square cavity whose each edge has the length $l=\np[m]{1.0}$. Body forces on the system are negligible. The constitutive elastic response of the disk is as follows:
\begin{equation}
\label{eq: INH1}
\tensor{P} = -p_{\s} \tensor{I} + G^{e} \tensor{F}. 
\end{equation}
We have used the following parameters: $\rho=\np[kg/m^{3}]{1.0}$, dynamic viscosities $\mu_{\f} = \mu_{\s} = \mu=\np[Pa\!\cdot\! s]{0.01}$, elastic shear modulus $G^{e} = \np[Pa]{0.1}$ and $U=\np[m/s]{1.0}$. For our numerical simulations we have used $\mathcal{Q}_{0}^{2}$ elements to represent $\bv{w}$ of the disk whereas we have used $\mathcal{Q}_{0}^{2}\vert\mathcal{P}_{-1}^{1}$ element for the fluid. The disk is represented using \np{320}~cells with \np{2626}~DoFs and the control volume has \np{4096}~cells and \np{45570}~DoFs. The time step size $h=\np[s]{1e-2}$. We consider the time interval $0<t \leq \np[s]{8}$ during which the disk is lifted from its initial position along the left vertical wall, drawn along underneath the lid and finally dragged downwards along the right vertical wall of the cavity (see Fig.~\ref{fig:LDCFlowBall-DGP-ResStress}). As the disk trails beneath the lid, it experiences large shearing deformations (see Fig.~\ref{fig:LDCFlowBall-DGP-ResStress-Deformation}). Due to incompressibility, the disk should have retained its original area over the course of time. However, as shown in Fig.~\ref{fig:LDCFlowBall-DGP-ResStress-AreaChange}(a), from our numerical scheme we obtain an area change of the disk of about $4\%$.  As discussed in great detail by \cite{Griffith-2012-a}, immersed methods are prone to poor volume conservation and a number of approaches have been proposed in the literature to address this issue. The error shown in Fig.~\ref{fig:LDCFlowBall-DGP-ResStress-AreaChange}(a) is more than that in \cite{GriffithLuo-2012-a} which is only of about 0.5\%. Note that  \cite{GriffithLuo-2012-a} use a finite difference scheme for solving the Navier-Stokes equation. The error in our scheme is certainly lower than the error ($> 20\%$) reported by \cite{WangZhang_2010_Interpolation-functions-0}, and on a par with the error ($\sim 5\%$) for the volume-conserving scheme used therein.  With this in mind, the volume conservation error in our method can be controlled with refinement as shown in Fig.~\ref{fig:LDCFlowBall-DGP-ResStress-AreaChange}(b), where the curve labeled `Case~2' is the repetition of that in Fig.~\ref{fig:LDCFlowBall-DGP-ResStress-AreaChange}(a) limited to the time interval where the error is largest, whereas the curves labeled `Case~1' and `Case~3' correspond to one level of refinement less and one higher, respectively, relative to the discretization used in Case~2.  The fact that the error in our method decreases with the increase in the refinement of the solid and the fluid domains is also demonstrated in Fig.~\ref{fig:pts-comp-solid-error-estimate} on p.~\pageref{fig:pts-comp-solid-error-estimate}, this time for the compressible case.
\begin{figure}[htbp]
	\begin{center}
		\includegraphics{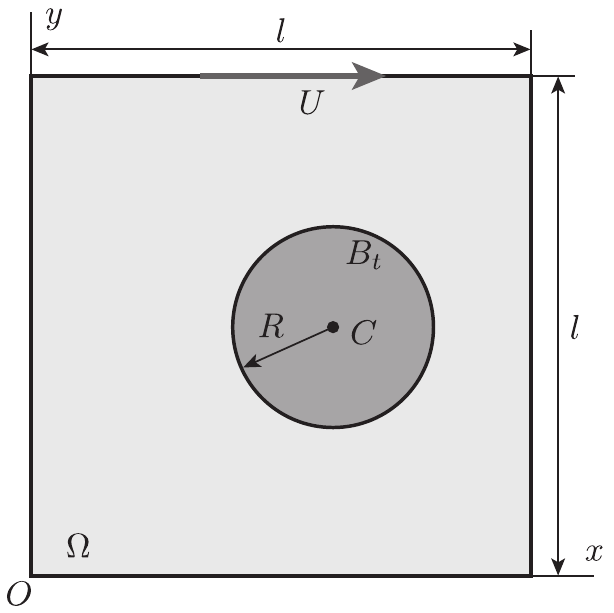}
		\caption{The initial configuration of an immersed disk entrained in a flow in a square cavity whose lid is driven with a velocity $U$ towards the right.}
	\label{fig:LDCFlow-Ball-Geometry}
	\end{center}
\end{figure}
\begin{figure}[htbp]
	\begin{center}
	\subfigure
	{\includegraphics[trim=15 10 10 20, clip=true]{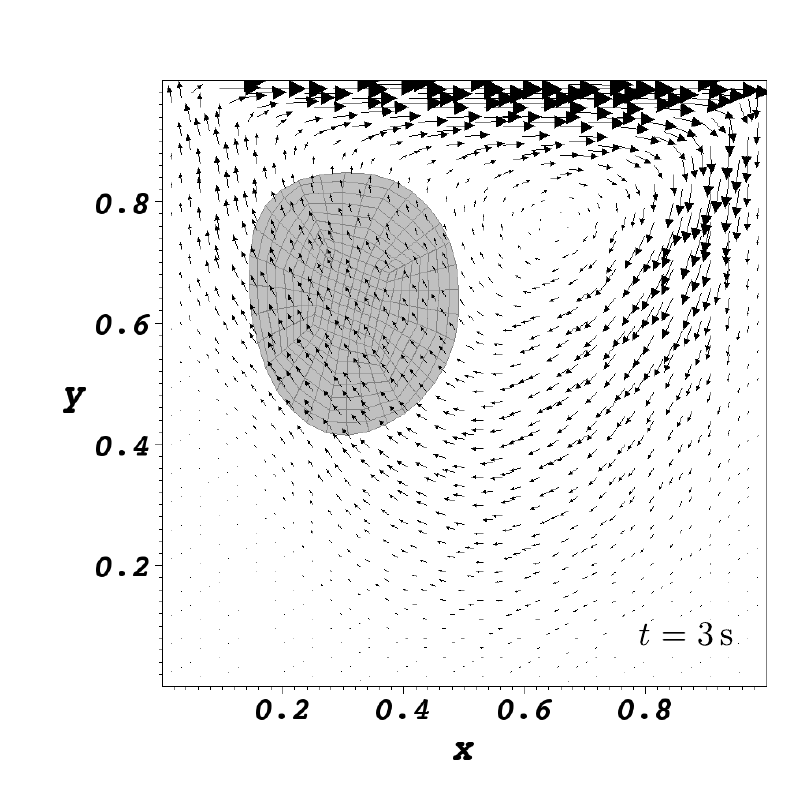}
	\label{fig:INH1-velocity-t3s}
	}
	\subfigure
	{\includegraphics[trim=15 10 10 20, clip=true]{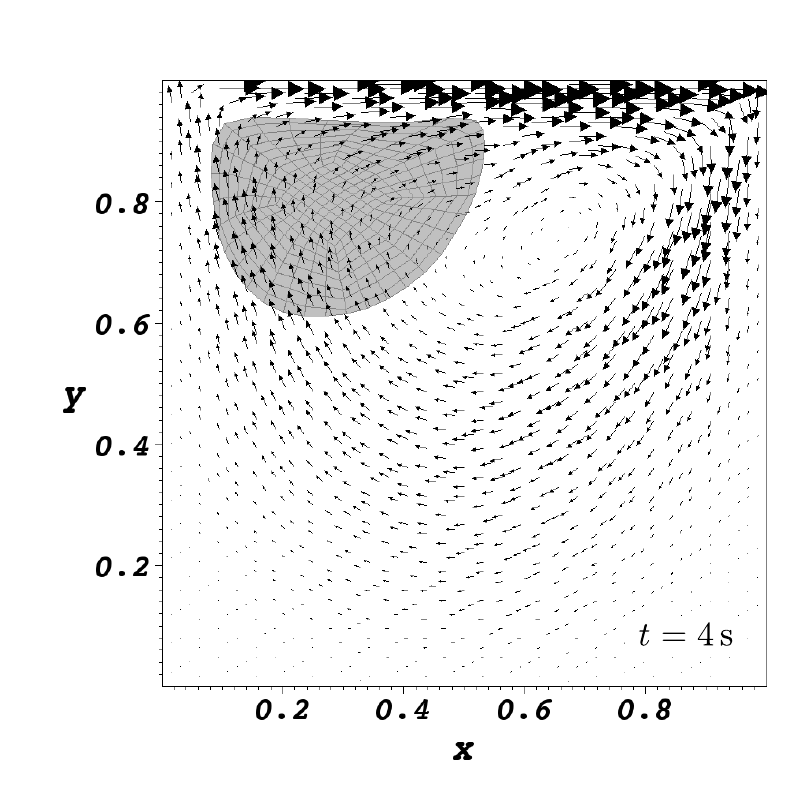}
	\label{fig:INH1-velocity-t4s}
	} 
	\subfigure
	{\includegraphics[trim=15 10 10 20, clip=true]{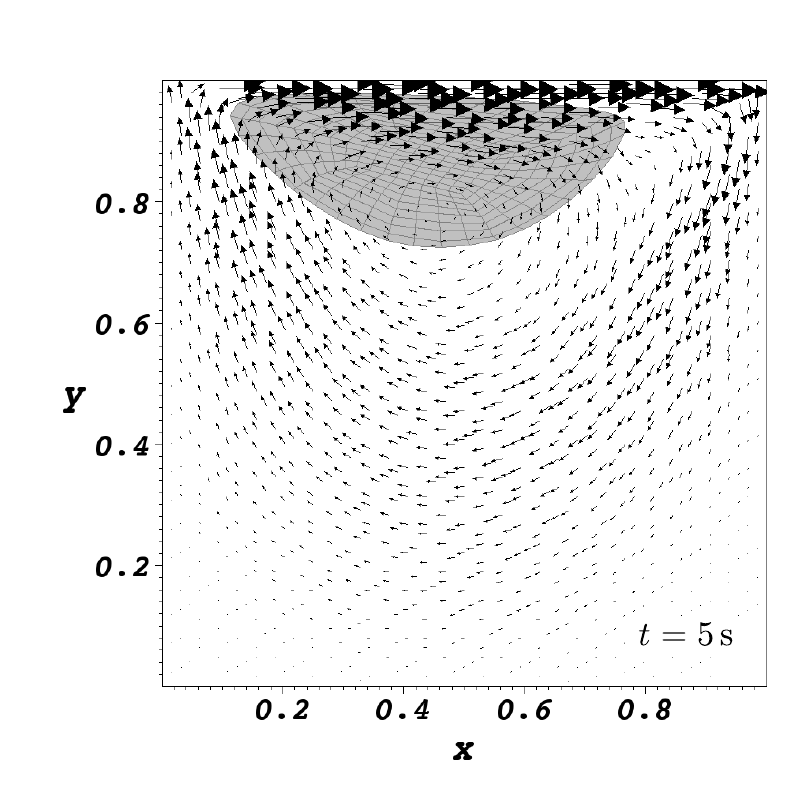}
	\label{fig:INH1-velocity-t5s}
	}
	\subfigure
	{\includegraphics[trim=15 10 10 20, clip=true]{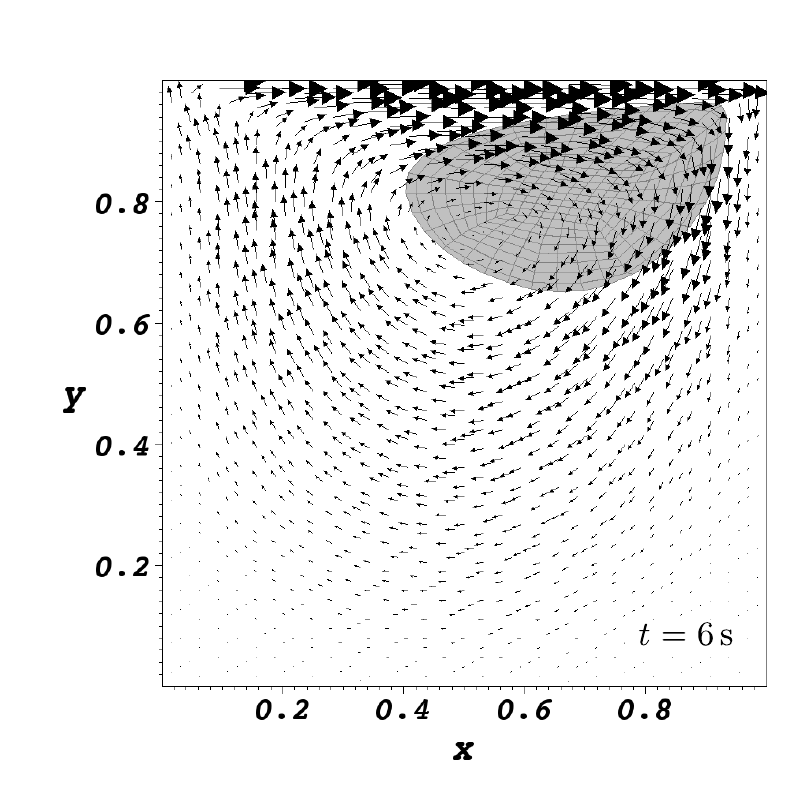}
	\label{fig:INH1-velocity-t6s}
	}
	\subfigure{
	\includegraphics[trim=15 10 10 20, clip=true]{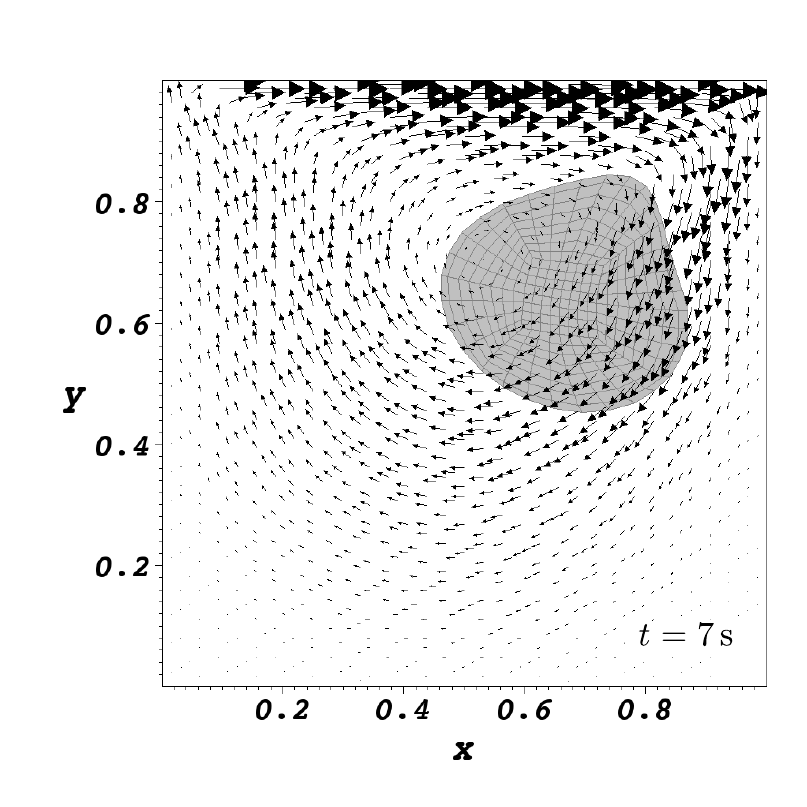}
	\label{fig:INH1-velocity-t7s}
	} 
	\subfigure{
	\includegraphics[trim=15 10 10 20, clip=true]{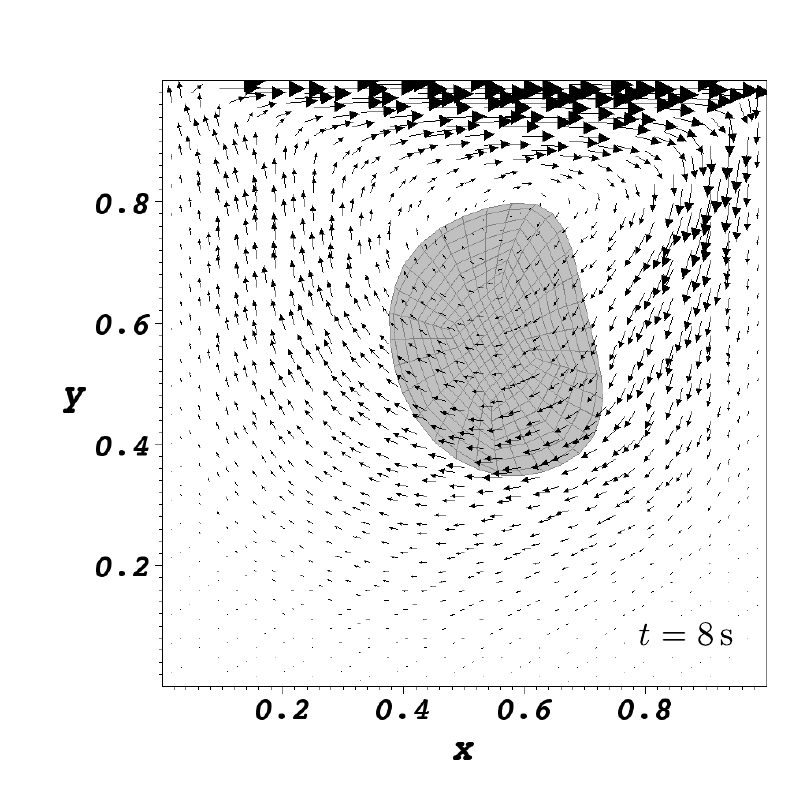}
	\label{fig:INH1-velocity-t8s}
	} 
\caption{The motion of a disk at different instants of time.}
\label{fig:LDCFlowBall-DGP-ResStress}
\end{center}
\end{figure}
\begin{figure}[htbp]
	\begin{center}
	\subfigure
	{\includegraphics[trim=0 10 0 70, clip=true]{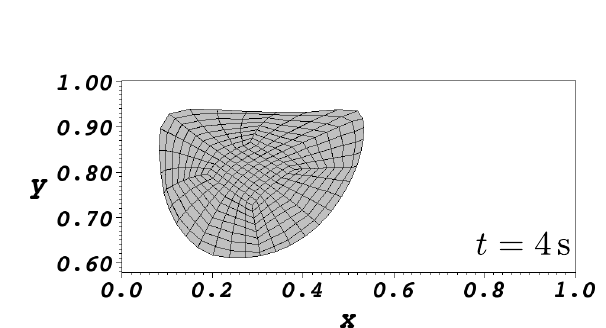} 
	\label{fig:INH1-solidlocation-t4s}
	}
	\subfigure
	{\includegraphics[trim=0 10 0 70, clip=true]{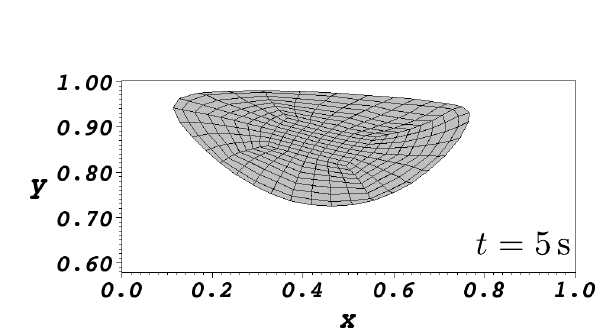}
	\label{fig:fig:INH1-solidlocation-t5s}
	} 
	\subfigure
	{\includegraphics[trim=0 10 0 70, clip=true]{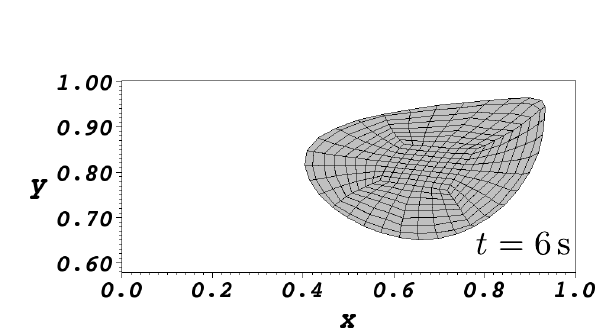}
	\label{fig:fig:INH1-solidlocation-t6s}
	}
	\caption{Enlarged view of the disk depicting its shape and location at various instants of time.}
	\label{fig:LDCFlowBall-DGP-ResStress-Deformation}
	\end{center}
\end{figure}
\begin{figure}[htbp]
\begin{center}
\includegraphics{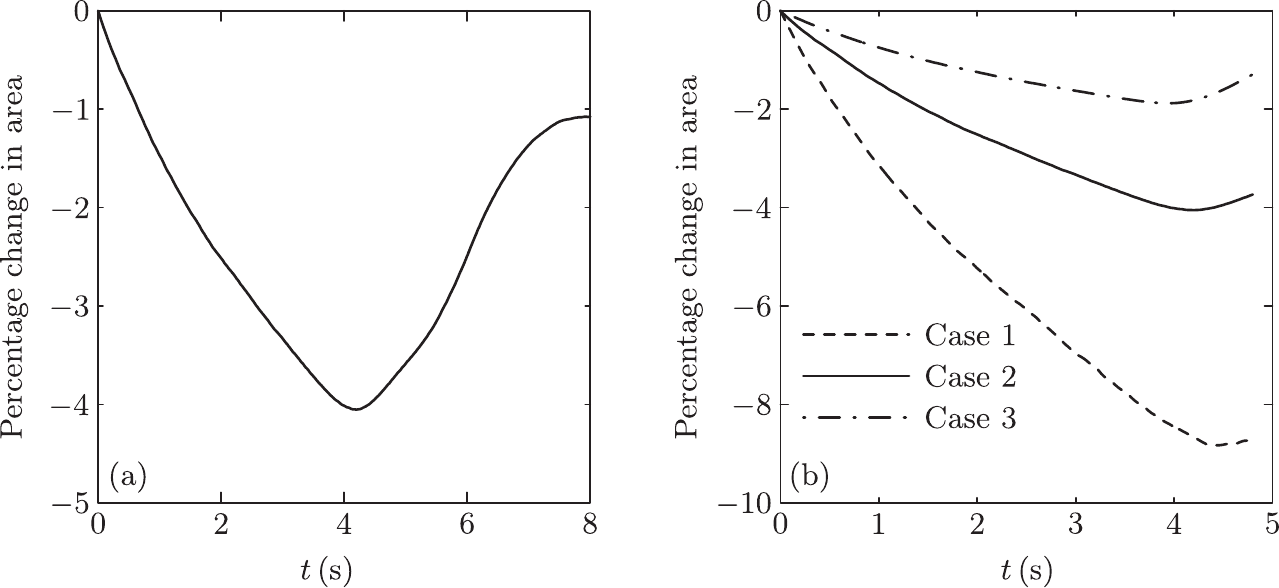}
\caption{(a) Percentage change in the area of the disk over time for the results shown in Fig.~\ref{fig:LDCFlowBall-DGP-ResStress} and~\ref{fig:LDCFlowBall-DGP-ResStress-Deformation}.  (b) Detail of plot (a) for three different refinement levels: Case 2 is the line in (a), whereas Cases 1 and 3 have one level of refinement lower and higher than that in Case 2, respectively.}
\label{fig:LDCFlowBall-DGP-ResStress-AreaChange}
\end{center}
\end{figure}

\subsubsection{Cylinder Falling in Viscous Fluid}
In the previous examples, the fluid and the solid had the same density and dynamic viscosity.  Furthermore, we consider the case of an incompressible solid with purely elastic constitutive response, i.e., without a viscous component in its stress constitutive law.  We present numerical tests partly inspired by those in \cite{ZhangGay-2007-Immersed-finite-0,ZhangGerstenberger-2004-Immersed-finite-0} meant to simulate experiments performed using a falling sphere viscometer. Specifically, we investigate the effect of changing the density of the solid and the viscosity of the fluid on the terminal velocity of the solid. Also, we estimate the drag coefficient of the falling object and discuss how the results depend on the mesh refinement.

The modeling of a cylinder sinking in a viscous fluid was first developed for the case of a rigid cylinder of length $L$ and radius $R$ released from rest in a quiescent unbounded fluid.  We use the suffix $\infty$ to distinguish the results for this ideal case from those concerning a cylinder sinking in a bounded fluid medium.  The longitudinal axis of the cylinder is assumed to be perpendicular to the direction of gravity and to remain so.  If the density of the cylinder and the fluid are the same, the cylinder will remain neutrally buoyant. If $\rho_{\s} > \rho_{\f}$, the weight of the cylinder will exceed the buoyancy force and will descend through the fluid with a velocity $u_{Cy\infty}$, parallel to the gravitational field, under the action of a net force with magnitude ${F}_{W}=\left( \rho_{\s} - \rho_{\f}\right) A L g$, where $A =\pi R^{2}$ is the cross-sectional area of the cylinder and $g$ is the acceleration due to gravity. Under creeping flow conditions, ${F}_{D\infty}= - \mu_{\f} \, c \, {u}_{Cy\infty}$, where $c$ is a constant. For a long cylinder, $c = 4 \pi L/[\ln (2 E)+ 1 -\kappa]$ \citep{CliftGrace-1978-Bubbles-drops-0} where $E= L/(2 R)$ is the aspect ratio for the cylinder and $\kappa$ is a constant whose values have been determined to be $0.72$ and $0.80685$ by different researchers.  When $F_{D\infty}$ balances $F_{W}$, the cylinder reaches terminal velocity ${U}_{t\infty} = [(\rho_{\s} - \rho_{\f}) A L g]/(\mu_{\f} \, c)$.  In real experiments, the fluid volume is finite and the resulting drag force $F_{D}$ is larger than $F_{D\infty}$ due to effect of the confining walls. With this in mind, following \cite{Brenner-1962-Effect-of-finite-0}, we have that $F_{D\infty}$ at some speed $u_{\infty}$ is related to $F_{D}$ at the same speed $u_{\infty}$ as $F_{D} = F_{D\infty}/\alpha$, where $0 < \alpha <1$ is a constant for a given experimental setup. Moreover, \cite{JayaweeraMason-1965-The-behaviour-of-freely-0} show that the terminal velocity in an unbounded medium $U_{t\infty}$ is related to its counterpart in a confined medium $U_{t}$ as follows: $U_{t}/U_{t\infty} = \alpha$, which implies that $U_{t} \propto (\rho_{\s} - \rho_{\f})/\mu_{\f}$.  This is the proportionality relation we expect to find in our numerical experiments.

We will express our results in terms of the flow's Reynolds number $\reynolds$ and the drag coefficient $C_{D}$, which, for the problem at hand, are defined, respectively, as
\begin{equation} \label{eqn:cylinder-Re-CD-bounded}
\reynolds = (2 R) \rho_{\f} u_{Cy} /\mu_{\f}
\qquad \mbox{and} \qquad
C_{D}= F_{D}\Big/\left[\tfrac{1}{2} \rho_{\f}\, u_{Cy}^{2} (2 R) L \right],
\end{equation}
which can be shown to imply 
$C_{D} \propto 1/\reynolds$
(see, \citealp{JayaweeraMason-1965-The-behaviour-of-freely-0}).


Referring to Fig.~\ref{fig:StokesFlow-Ball-Geometry},
\begin{figure}[htbp]
	\centering
	  \includegraphics[scale=0.9]{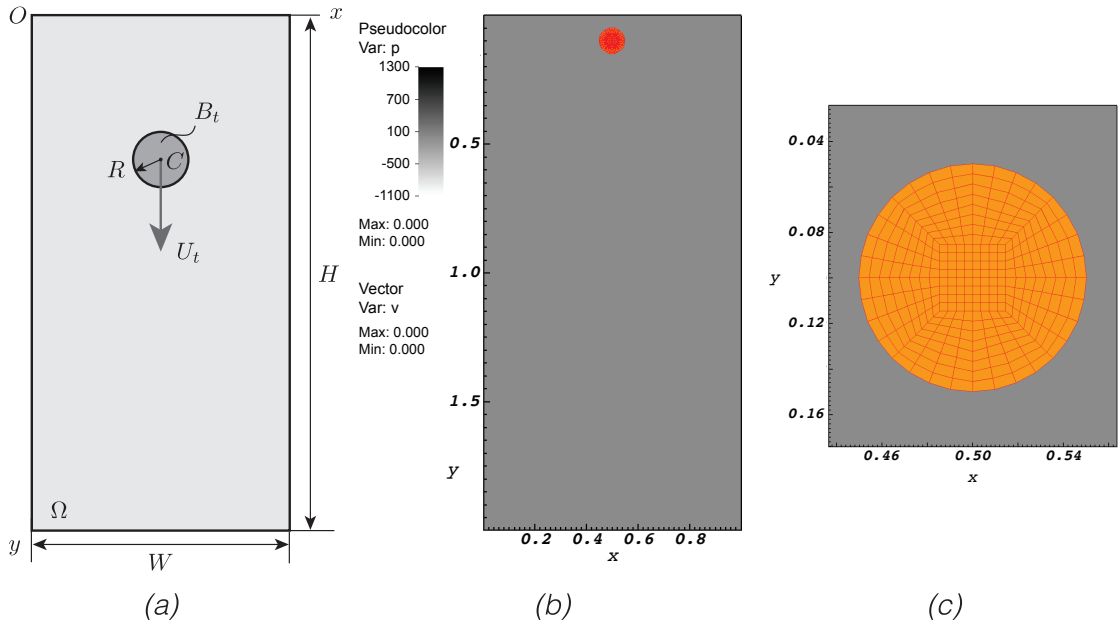}
	\caption{Vertical channel containing viscous fluid through which a small incompressible disk is descending: (a) system's geometry; (b) initial conditions (disk released from rest in a quiescent fluid); (c) detail of mesh of the immersed body. The disk's terminal velocity is denoted by $U_{t}$.}
	\label{fig:StokesFlow-Ball-Geometry}
\end{figure}
in our numerical experiments we consider a disk $B_{t}$ (representing the midplane of the cylinder), with radius $R$ and center $C$, released from rest in an initially quiescent rectangular control volume $\Omega$ with height $H$ and width $W$.  As the disk sinks, we measure the position of $C$ and infer its velocity, denoted by $u_{Cy}$.  When $u_{Cy}$ achieves a (sufficiently) constant value, we refer to this value as the terminal velocity and denote it by $U_{tN}$.  When $u_{Cy} = U_{tN}$ we also compute the drag coefficient $C_{DN}$ of the disk.  The latter is assumed to consists of an incompressible neo-Hookean material whose Piola stress is given by
\begin{equation} 
\label{eq: INH0}
\tensor{P}_{\s} = -p_{\s} \tensor{I} + G^{e} \bigl(\tensor{F} -\invtrans{\tensor{F}}\bigr).
\end{equation}
The following parameters were used: $H = \np[cm]{2.0}$, $W = \np[cm]{1.0}$, $R=\np[cm]{0.05}$, $\rho_{\f}=\np[g/cm^{3}]{1.0}$, $\mu_{\s}=0$, $G^{e} = \np[dyn/cm^{2}]{1e3}$ and $g=\np[cm/s^{2}]{981}$. Two different cases were considered. In Case~1, the density of the solid was varied while the viscosity of the fluid was held constant at $\mu_{\f}=\np[P]{1.0}$. In Case~2, we used for different fluid viscosities while the density of the solid was held constant at $\rho_{\s}= \np[g/cm^{3}]{3.0}$. The values used for $\rho_{\s}$ and $\mu_{\f}$ can be found in Table~\ref{tab:properties-falling-cylinder},
\begin{table}[htbp]\small
\caption{Values of solid density and fluid density used for the tests. Also shown are the values obtained for the terminal velocity, Reynolds number and drag coefficient from these tests.}
\begin{center}
\begin{tabular*} {0.95\textwidth} {@{\extracolsep{\fill}} c c c c c c}
\toprule
&
$\rho_{\s} (\text{g}/\text{cm}^{3})$ &
$\mu_{\f} (\text{P})$  &
$U_{tN} (\text{cm}/\text{s})$ & 
$\reynolds_{tN}$ &
$C_{DtN}$
\\
\midrule
\multirow{3}{*}{Case 1}
&	 2.0 &	 1.0 &	 0.8179	& 0.082 & 230.3\\
&	 3.0 &	 1.0 &	 1.6270 & 0.163 & 116.4\\
& 	 4.0 &	 1.0 &	 2.4412 & 0.244 &  77.6\\
\midrule
\multirow{4}{*}{Case 2}
& 3.0 & 1.0 & 1.6270 & 0.163 & 	 116.4 \\
& 3.0 & 2.0 & 0.8236 & 0.041 & 	 454.4\\
& 3.0 & 4.0 & 0.4137 & 0.010 & 1801.0\\
& 3.0 & 8.0 & 0.2070 & 0.003 & 7192.5\\
\bottomrule
\end{tabular*}
\end{center}
\label{tab:properties-falling-cylinder}
\end{table}
which also lists the corresponding values of $U_{tN}$, $\reynolds_{N}$, and $C_{DtN}$. We have used $\mathcal{Q}_{0}^{2}\vert\mathcal{Q}_{0}^{1}$ elements for the control volume and a Gauss quadrature rule of order $4$ for assembling the operators defined over the solid domain. In all numerical experiments pertaining to Cases~1 and~2, we have used \np{16384}~cells and \np{181250}~DoFs for the control volume, and \np{320}~cells and \np{2626}~DoFs for the solid.

The results for Case~1 are reported in Figs.~\ref{fig:Density-change1} and~\ref{fig:Density-change2},
\begin{figure}[htbp]
	\centering
	\subfigure[Axial position of mass center of the disk.]
	{\includegraphics[height=2.17in]
	{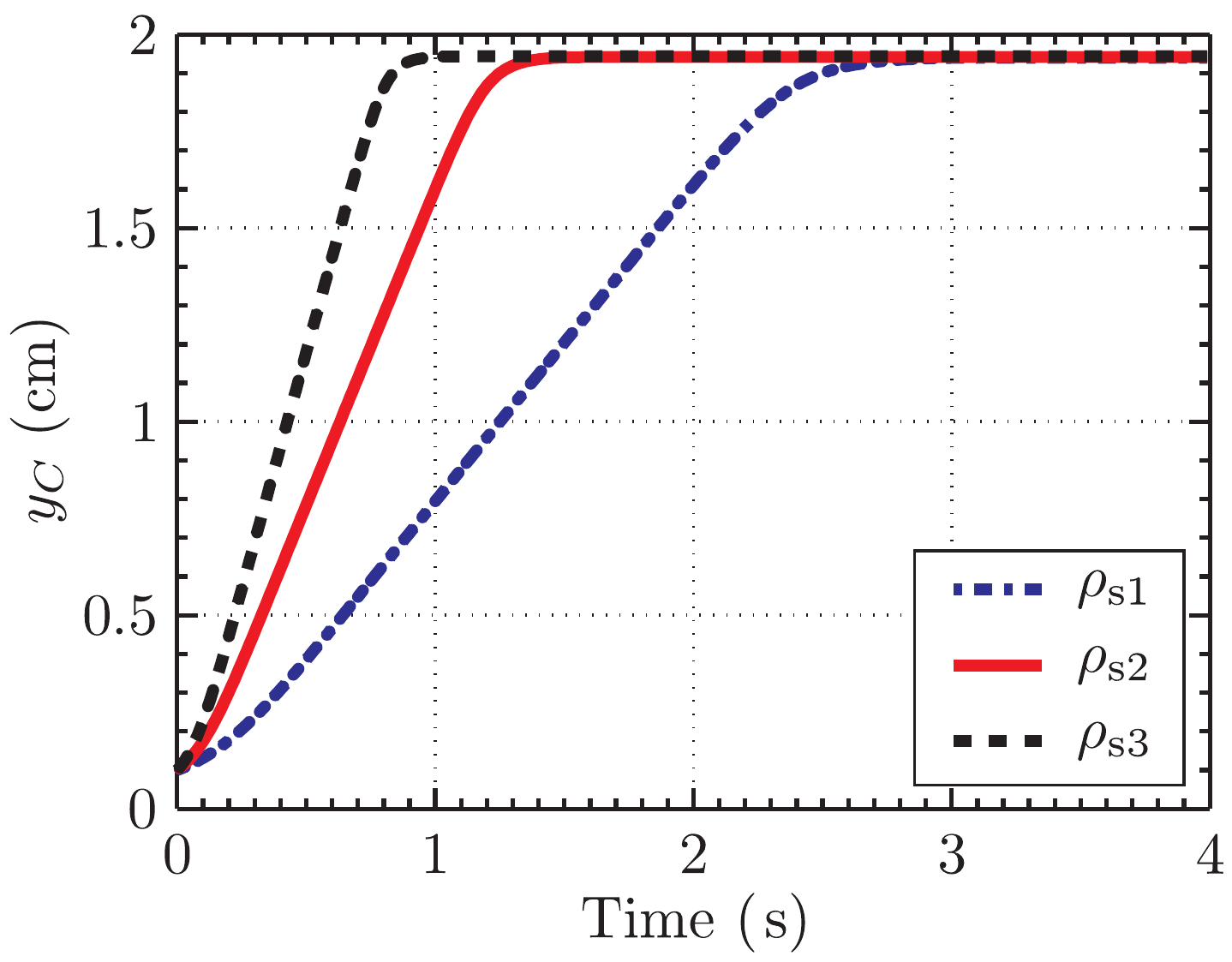}
	\label{subfig:TermVel-y-versus-t-RhosVarying}
	}
	\subfigure[Vertical velocity of mass center of the disk.]
	{\includegraphics[height=2.17in]
	{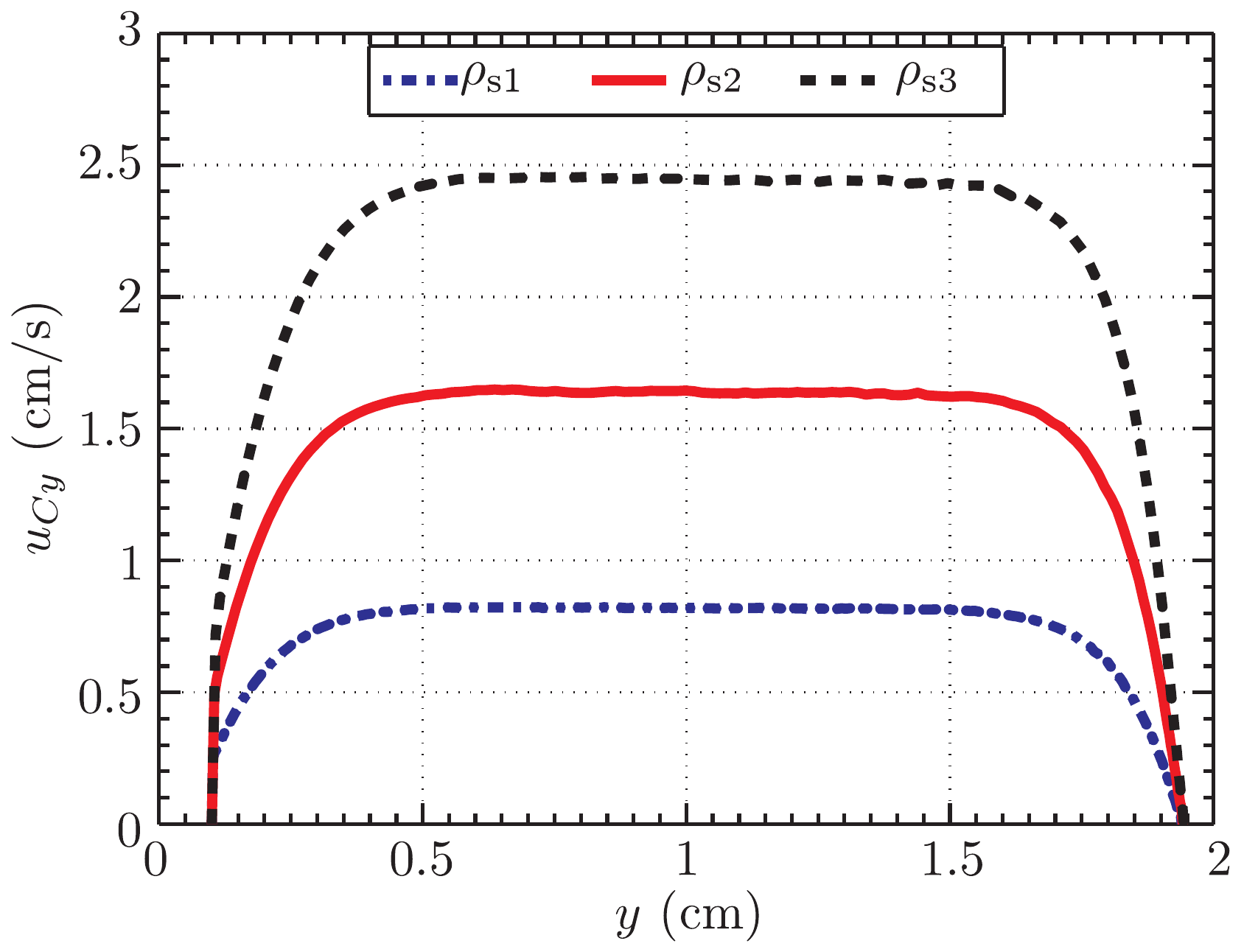}
	\label{subfig:TermVel-v-versus-y-RhosVarying}
	}
	\caption{Effect of changing the density of the solid on the motion of the center of mass of the cylinder. Note: $\rho_{\s1} = {\np[g/cm^{3}]{2}}$, $\rho_{\s2} = {\np[g/cm^{3}]{3}}$ and $\rho_{\s3} = {\np[g/cm^{3}]{4}}$.}
	\label{fig:Density-change1}
\end{figure}
\begin{figure}[htbp]
	\centering
	\subfigure[Velocity of mass center of the cylinder versus time.]
	{\includegraphics[width=2.82in]
	{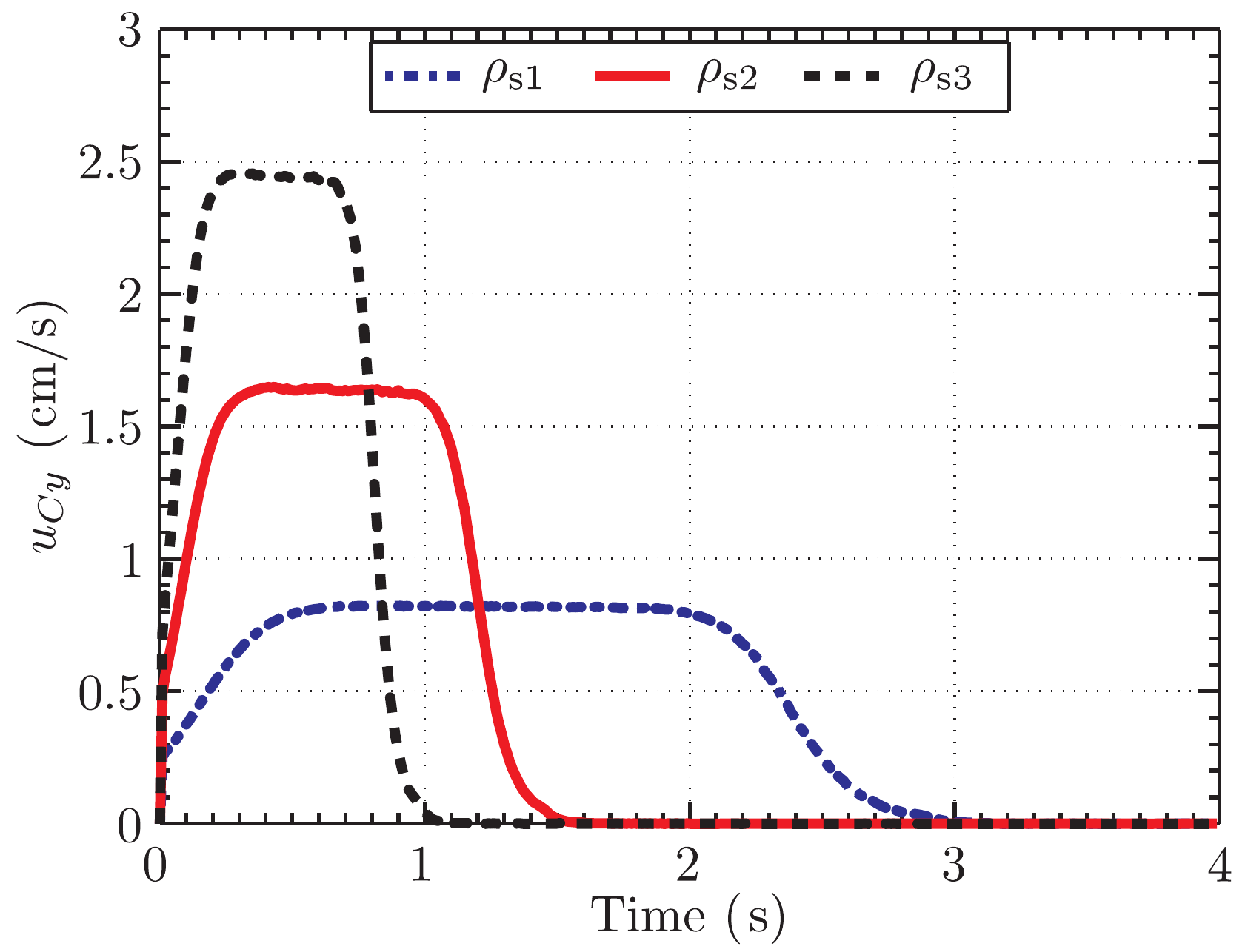}
	\label{subfig:TermVel-v-versus-t-RhosVarying}
	}
	\subfigure[Terminal velocity of the cylinder as function of its density]
	{\includegraphics[width=2.82in]
	{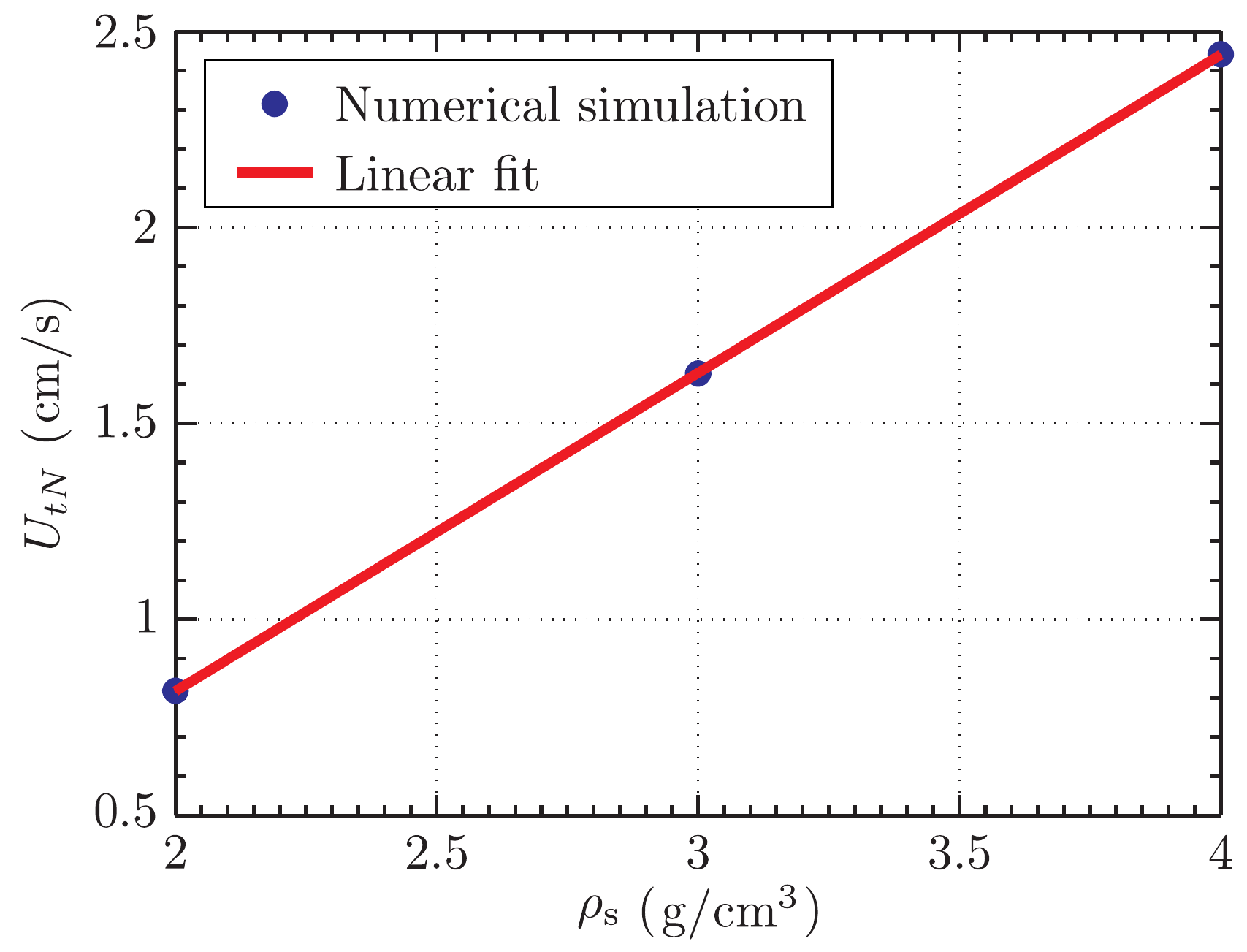}
	\label{subfig:TermVel-Effect-of-Density-on-TV}
	}
	\caption{Effect of changing the density of the cylinder on its terminal velocity. 	Note: $\rho_{\s1}={\np[g/cm^{3}]{2}}$, $\rho_{\s2}={\np[g/cm^{3}]{3}}$ and $\rho_{\s3}={\np[g/cm^{3}]{4}}$.
}
	\label{fig:Density-change2}
\end{figure}
whereas those for Case~2 are reported in Figs.~\ref{fig:Viscosity-change1} and~\ref{fig:Viscosity-change2}.
\begin{figure}[htbp]
	\centering
	\subfigure[Axial position of mass center of the disk.]
	{\includegraphics[width=2.82in]
	{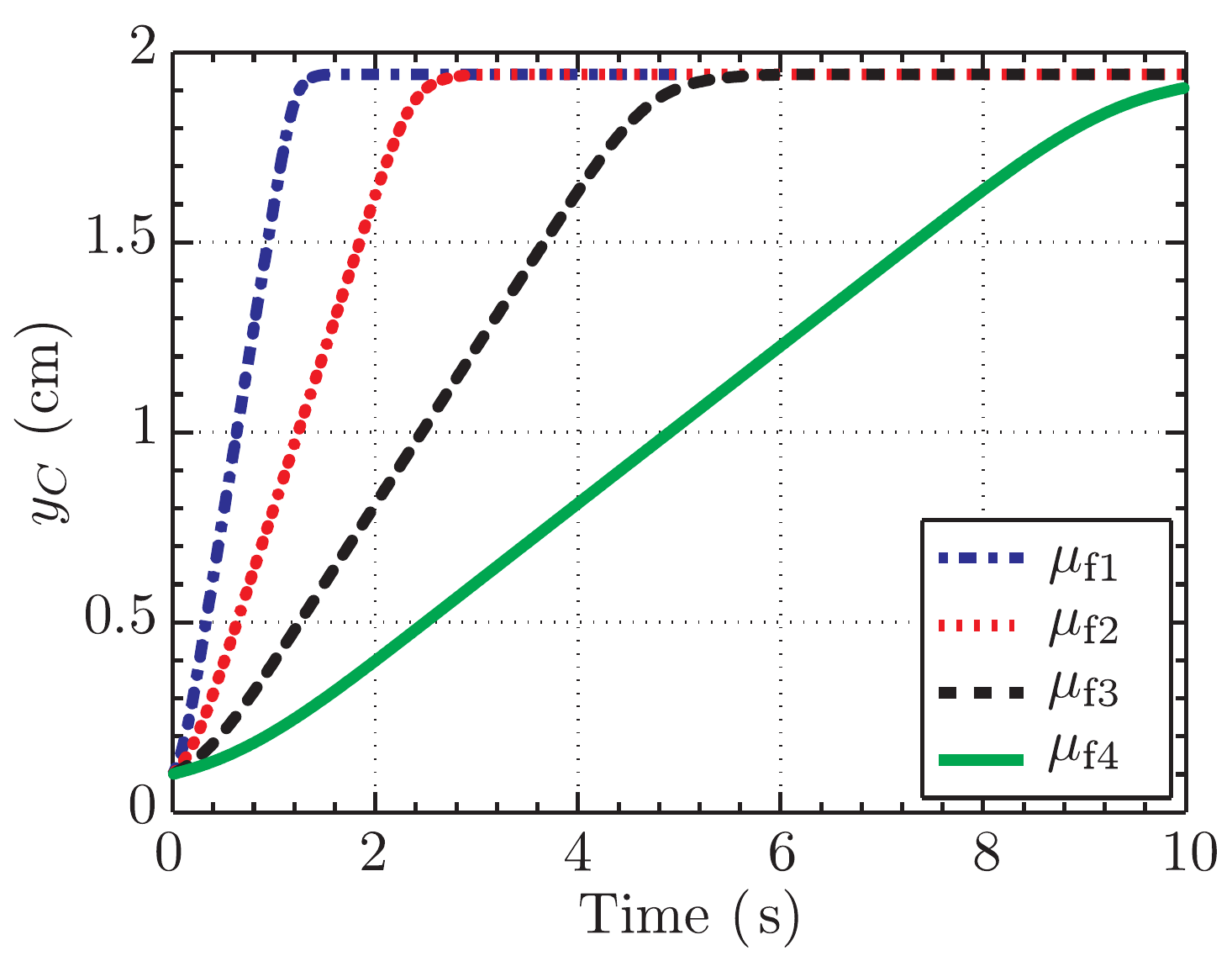}
	\label{subfig:TermVel-y-versus-t-EtafVarying}
	}
	\subfigure[Vertical velocity of mass center of the disk]
	{\includegraphics[width=2.82in]
	{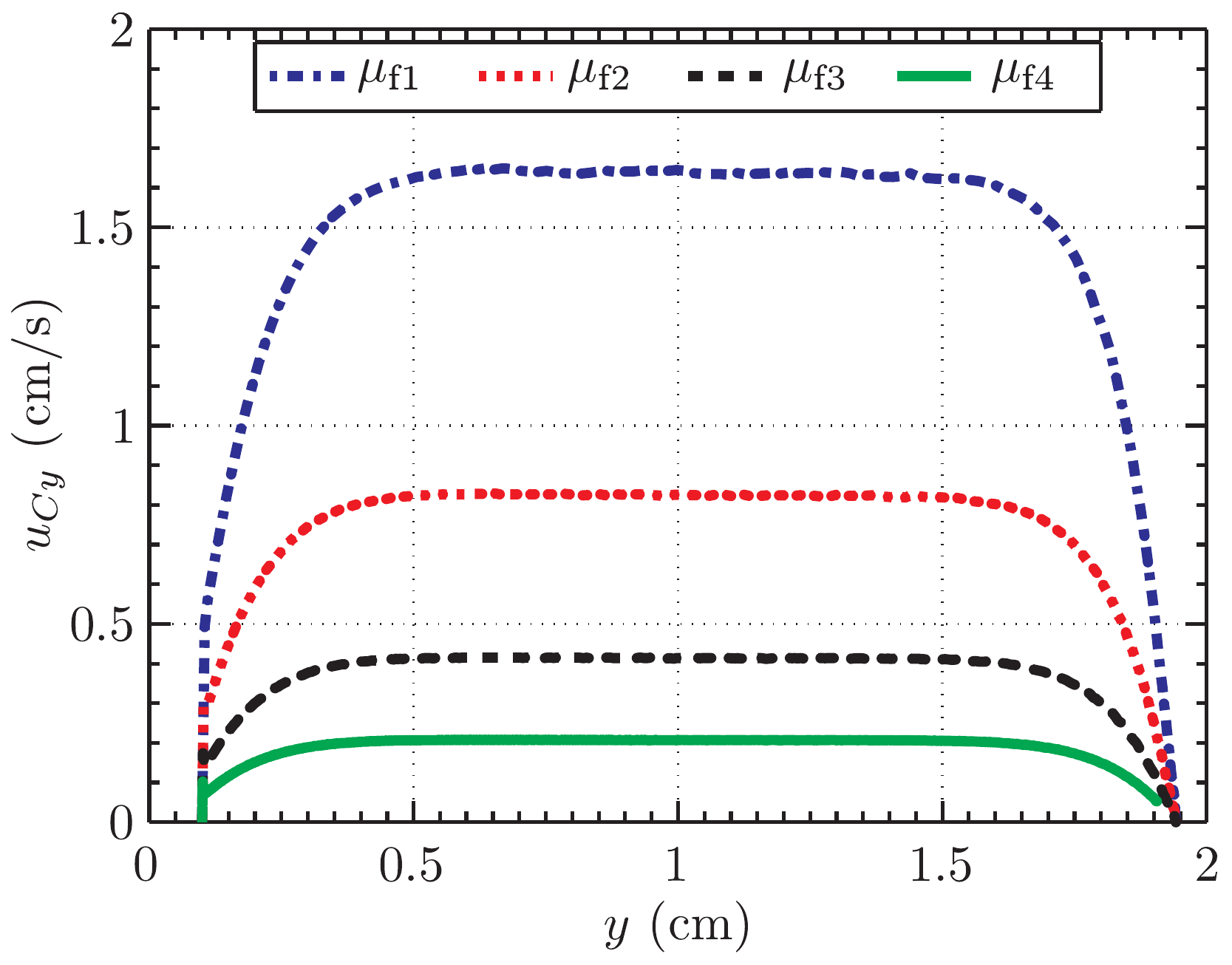}
	\label{subfig:TermVel-v-versus-y-EtafVarying}
	}
	\caption{Effect of changing the viscosity of the fluid on the motion of the center of mass of the cylinder. Note: $\mu_{\f1}=\np[P]{1.0}$, $\mu_{\f2}=\np[P]{2.0}$, $\mu_{\f2}=\np[P]{4.0}$ and $\mu_{\f4}=\np[P]{8.0}$.}
	\label{fig:Viscosity-change1}
\end{figure}
\begin{figure}[htbp]
	\centering
	\subfigure[Velocity of mass center of the cylinder versus time.]
	{
	\includegraphics[width=0.47\textwidth]
	{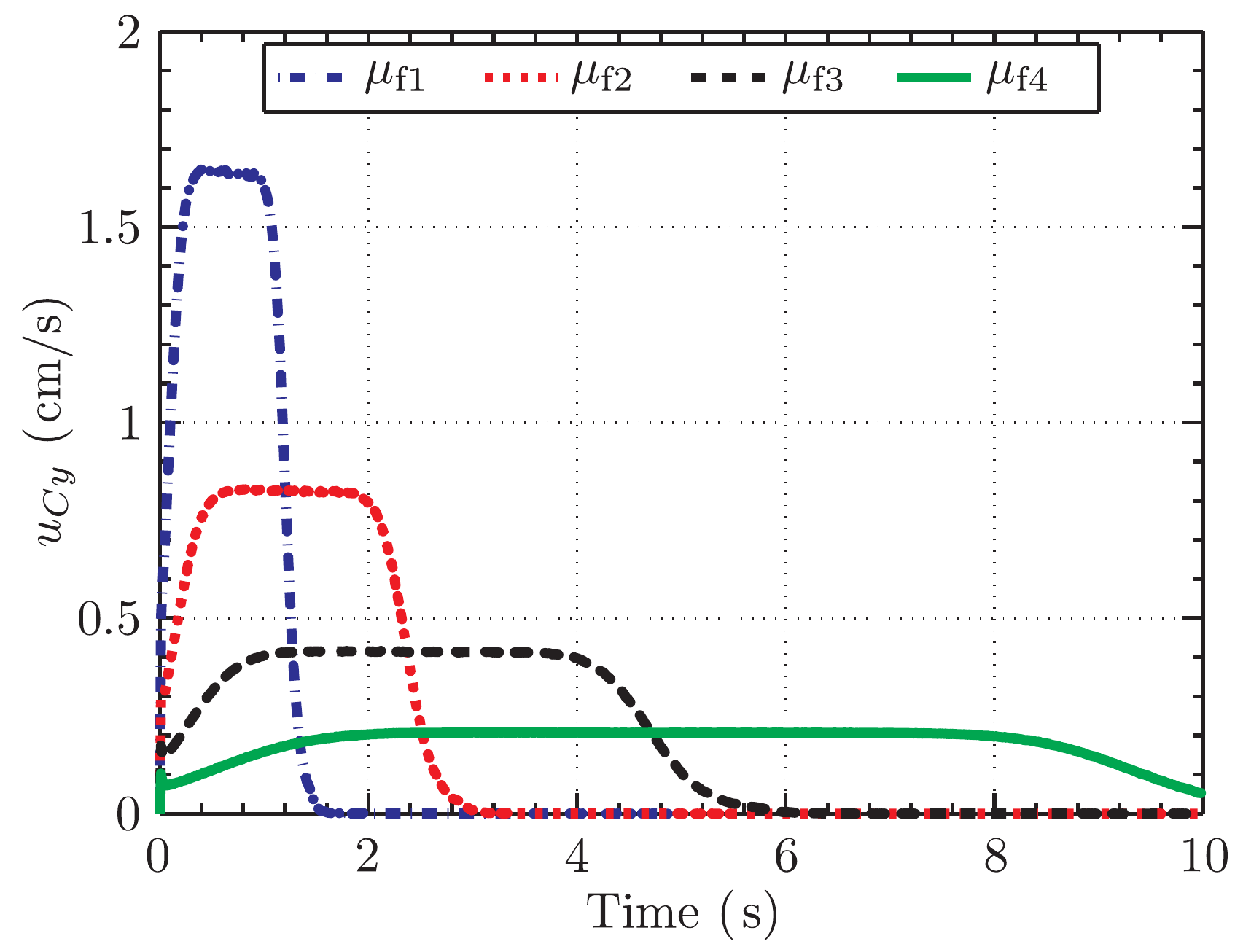}
	\label{subfig:TermVel-v-versus-t-EtafVarying}
	}
	\subfigure[Terminal velocity of the cylinder as function of the viscosity of the fluid.]
	{
	\includegraphics[width=0.47\textwidth]
	{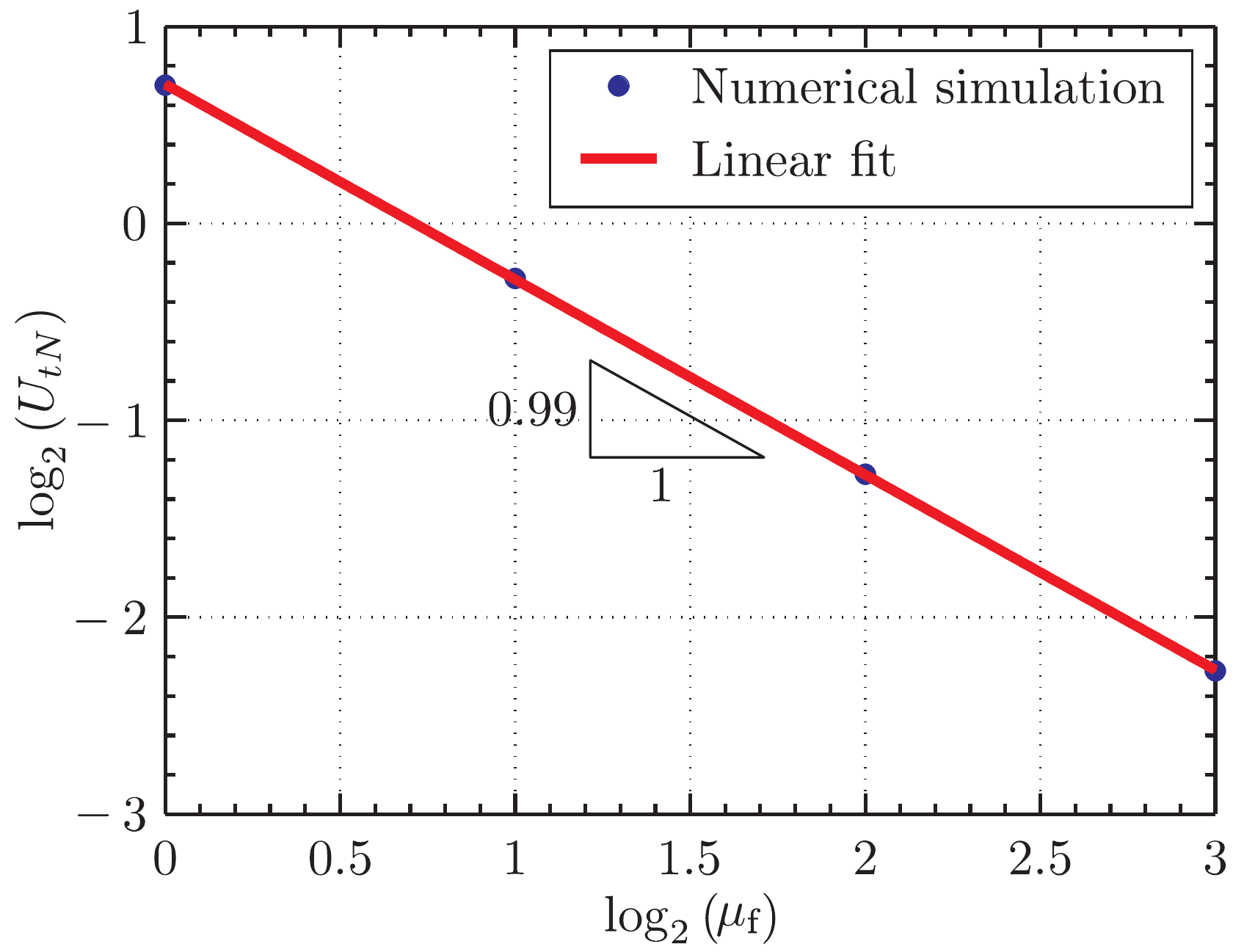}
	\label{subfig:TermVel-Effect-of-Etaf-on-TV}
	}
	{\caption{Effect of fluid viscosity on the terminal velocity of the disk. Note: $\mu_{\f1}={\np[P]{1.0}}$, $\mu_{\f2}={\np[P]{2.0}}$, $\mu_{\f2}={\np[P]{4.0}}$ and $\mu_{\f4}={\np[P]{8.0}}$.}}
	\label{fig:Viscosity-change2}
\end{figure}
Figures~\ref{subfig:TermVel-v-versus-y-RhosVarying} and \ref{subfig:TermVel-v-versus-y-EtafVarying} show that $u_{Cy}$ increases over a distance of about $y<H/3$, remains constant over $H/3<y<2H/3$ and then decreases over the remaining length of the channel. This is in accordance with actual experimental observations (cf., \citealp{CliftGrace-1978-Bubbles-drops-0}). From Figs.~\ref{subfig:TermVel-v-versus-t-RhosVarying} and~\ref{subfig:TermVel-v-versus-t-EtafVarying}, we see that the disk's terminal velocity increases with its density and decreases with the increase in the fluid's viscosity. Moreover, from Fig.~\ref{subfig:TermVel-Effect-of-Density-on-TV} we see that the terminal velocity is linearly proportional to the density of the disk, and from Fig.~\ref{subfig:TermVel-Effect-of-Etaf-on-TV} we see that the terminal velocity is inversely proportional to the viscosity of the fluid, as expected. Finally, from Fig.~\ref{fig:TermVel_Cd-vs-Re_fr7sr3}
\begin{figure}[htbp]
	\centering
	\includegraphics[width=3.5in]
	{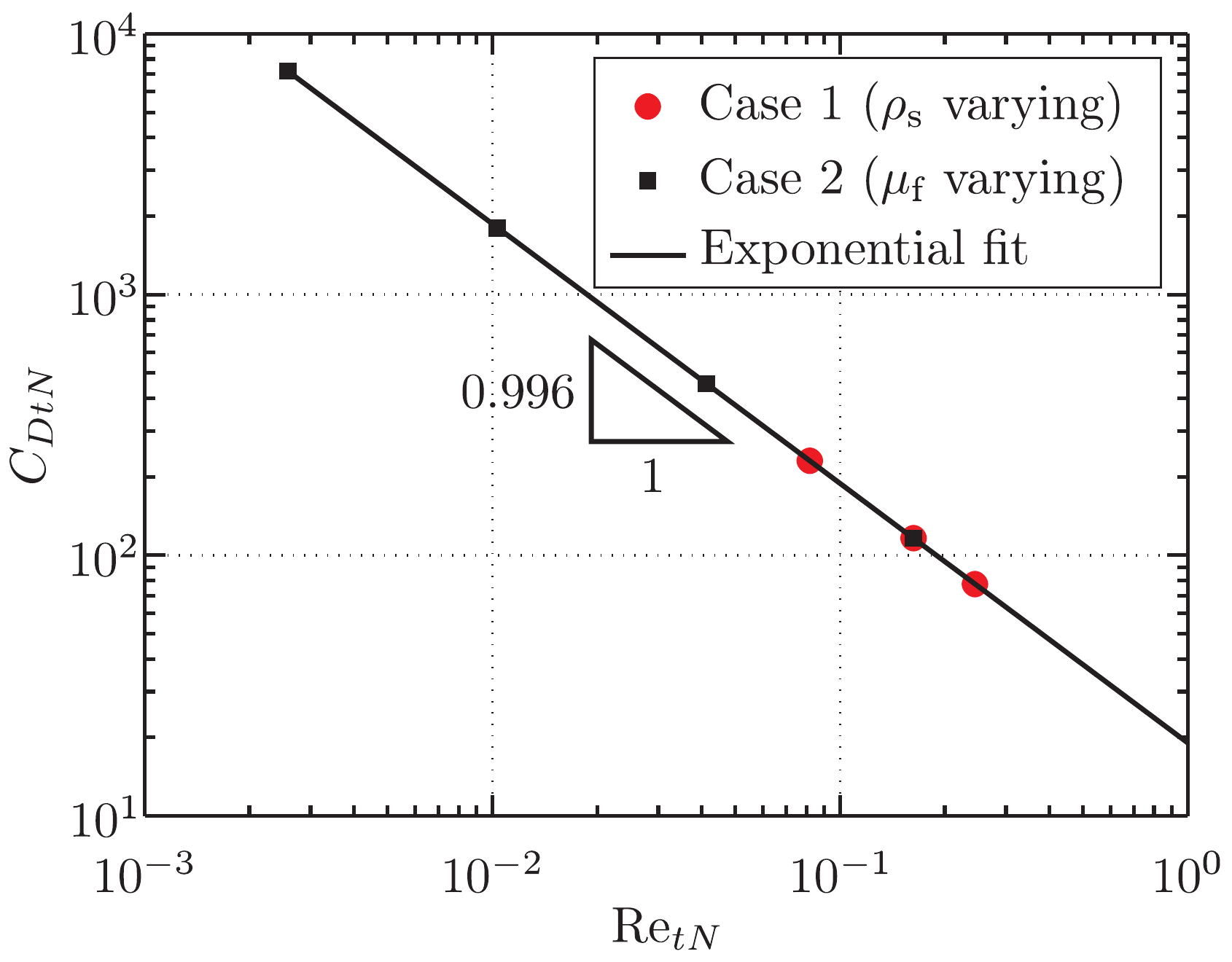}
	\caption{Drag coefficient versus Reynolds number (corresponding to the terminal velocity of the cylinder).}
	\label{fig:TermVel_Cd-vs-Re_fr7sr3}
\end{figure}
we see that the calculated drag coefficient at the terminal velocity $C_{DtN}$ of the cylinder is indeed inversely proportional to the corresponding Reynolds number $\reynolds_{tN}$.

We end this section with a few remarks on convergence and mesh refinement. The computed value of the terminal velocity can be expected to be accurate only when the meshes are ``sufficiently'' refined. With this in mind, we considered the effect of mesh refinement on the the value of $U_{tN}$ for Case~1 corresponding to $\rho_{\s}=\np[g/cm^{3}]{4}$. Table~\ref{tab:mesh-convergence-falling-cylinder} shows the mesh sizes for both the solid and the control volume along with the corresponding value of $U_{tN}$.
\begin{table}[htbp]\small
\caption{Terminal velocity of the cylinder obtained from simulations using meshes having different global refinement levels.}
\begin{center}
\begin{tabular*} {0.95\textwidth} {@{\extracolsep{\fill}} c c c c c c}
\toprule
&
\multicolumn{2}{c}{Solid}&
\multicolumn{2}{c}{Control Volume} &
\multirow{2}{*}{$U_{tN}{\text{cm}/\text{s}}$} 
\\
\cmidrule(r){2-3}\cmidrule(r){4-5}
&
Cells &
DoFs &
Cells &
DoFs &
\\
\midrule
Level 1&	   \np{80} &	    \np{674} &	 \np{1024} &	 \np{11522} & 2.10 \\
Level 2&	   \np{80} &    \np{674} &	 \np{4096} &  \np{45570} & 2.37 \\
Level 3&  \np{320} &   \np{2626} &	\np{16384} & \np{181250} & 2.45 \\
Level 4& \np{1280} &  \np{10370} & 	\np{65536} & \np{722946} & 2.50 \\
\bottomrule
\end{tabular*}
\end{center}
\label{tab:mesh-convergence-falling-cylinder}
\end{table}
The corresponding velocity of $C$ as a function of time is shown in Fig.~\ref{fig:TermVel_mesh-convergence},
\begin{figure}[htbp]
	\centering
	\includegraphics[width=3.5in]
	{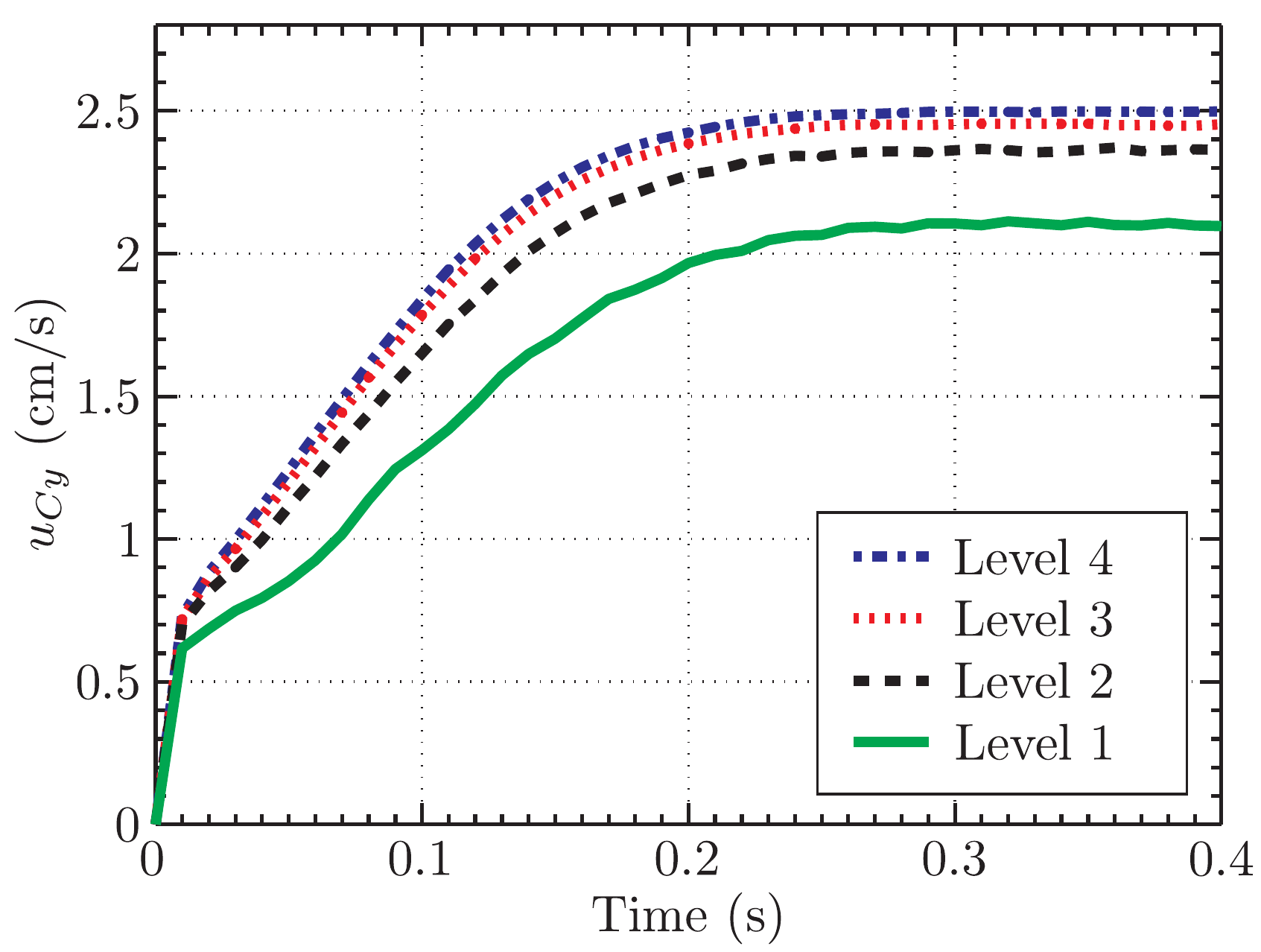}
	\caption{Effect of mesh refinement level on the terminal velocity of the mass center of the cylinder.}
	\label{fig:TermVel_mesh-convergence}
\end{figure}
in which we see that, as the meshes are refined, $U_{tN}$ tends to achieve a ``converged'' value. We note that the Case~1 and~2 results presented earlier correspond to Level~3 in Table~\ref{tab:mesh-convergence-falling-cylinder}.

\subsection{Results for Compressible Immersed Solids}
\subsubsection{Compressible Annulus inflated by a Point source}
Here we present a problem involving a compressible solid with stress response containing both elastic and viscous contributions. Specifically, we study the deformation of a hollow compressible cylinder submerged in a fluid contained in a rigid prismatic box due to the influx of fluid along the axis of the cylinder. Referring to Fig.~\ref{fig:PointSource-Annulus-Geometry},
\begin{figure}[htbp]
	\centering
		\includegraphics{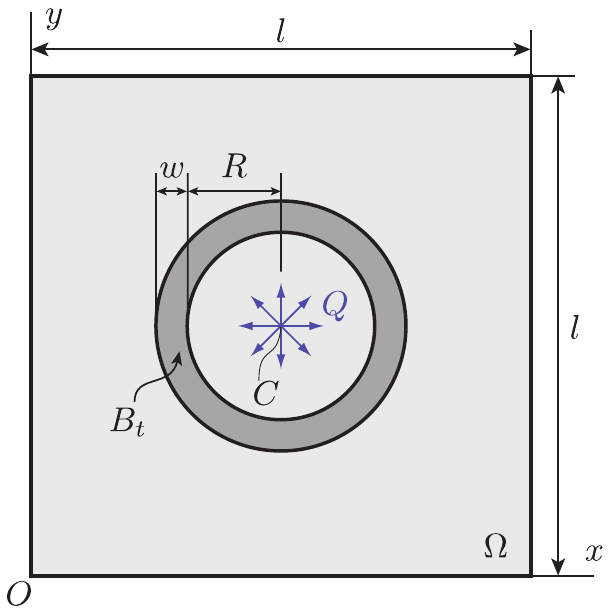}
		\caption{Initial configuration of an annulus immersed in a square box filled with fluid. At the center $C$ of the box is a point source of constant strength $Q$.}
		\label{fig:PointSource-Annulus-Geometry}
\end{figure}
we consider a two-dimensional solid annulus with inner radius $R$ and thickness $w$ that is concentric with a fluid-filled square box of edge length $l$. A point mass source of fluid of constant strength $Q$ is located at the center $C$ of $\Omega$.  Because of this source, the balance of mass for the system is modified as follows:
\begin{equation} \label{eqn:WeakBalOfMass-PtSource}
\int_{\Omega} q \biggl[ \nabla\cdot \bv{u} +
\frac{Q}{\rho_{\f}}\delta (\bv{x}-\bv{x}_{C}) \biggr] \d{v}
-
\int_{B_{t}} q \nabla\cdot \bv{u} \d{v}
= 0, 
\end{equation}
where $\delta \left(\bv{x}-\bv{x}_{C}\right)$ denotes a Dirac-$\delta$ distribution centered at $C$.  We apply homogeneous Dirichlet boundary conditions on $\partial \Omega$. The annulus was chosen to have a compressible Neo-Hookean elastic response given by 
\begin{equation} \label{eqn:FirstPKStress-CompSolid}
\tensor{P}^{e}_{\s}=G^{e} \bigl(\tensor{F}- J^{-2 \nu/(1-2\nu)} \tensor{F}^{-\mathrm{T}} \bigr),
\end{equation}
where $G^{e}$ is the elastic shear modulus and $\nu$ is the Poisson's ratio for the solid.  Since the solid is compressible both the volume of solid and that of the fluid in the control volume can change.   However, since the fluid cannot leave the control volume, the amount of fluid volume increase must match the decrease in the volume of the solid. This implies that the difference in these two volumes can serve as an estimate of the numerical error incurred.

We have used the following parameters: $R =\np[m]{0.25}$, $w=\np[m]{0.05}$, $l=\np[m]{1.0}$, $\rho_{\f}= \rho_{\s_{0}}= \np[kg/m^{3}]{1}$, $\mu_{\f}= \mu_{\s}= \np[Pa \! \cdot \! s]{1}$, $G^{e}=\np[Pa]{1}$, $\nu=0.3$, $Q=\np[kg/s]{0.1}$ and $dt=\np[s]{0.01}$.  We have tested three different mesh refinement levels whose details have been listed in Table~\ref{tab:mesh-refinement-pts-comp-solid}.
\begin{table}[htbp]\small
\caption{Number of cells and DoFs used in the different simulations involving the deformation of a compressible annulus under the action of a point source.}
\begin{center}
\begin{tabular*} {0.95\textwidth} {@{\extracolsep{\fill}} c c c c c}
\toprule
&
\multicolumn{2}{c}{Solid}&
\multicolumn{2}{c}{Control Volume} 
\\
\cmidrule(r){2-3}\cmidrule(r){4-5}
&
Cells &
DoFs &
Cells &
DoFs
\\
\midrule
Level 1&  \np{6240} &   \np{50960} &  \np{1024} &   \np{9539} \\
Level 2& \np{24960} &  \np{201760} &  \np{4096} &  \np{37507} \\
Level 3& \np{99840} &  \np{802880} & \np{16384} & \np{148739} \\
\bottomrule
\end{tabular*}
\end{center}
\label{tab:mesh-refinement-pts-comp-solid}
\end{table}
The initial state of the system is shown in Fig.~\ref{subfig:pts-comp-solid-t0}.
\begin{figure}[htbp]
\centering
\subfigure[$t={\np[s]{0}}$]
{\includegraphics[width=0.48\textwidth]{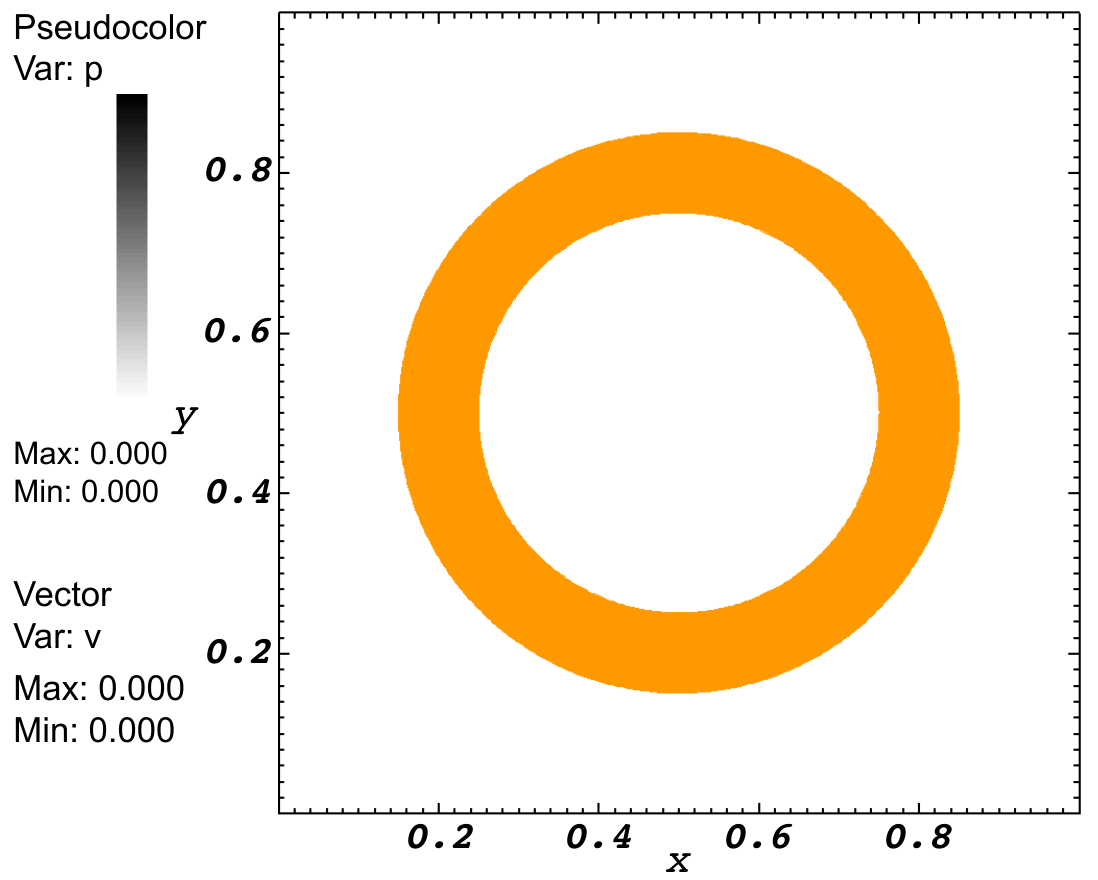}
\label{subfig:pts-comp-solid-t0}
}
\subfigure[$t={\np[s]{1.0}}$]
{\includegraphics[width=0.48\textwidth]{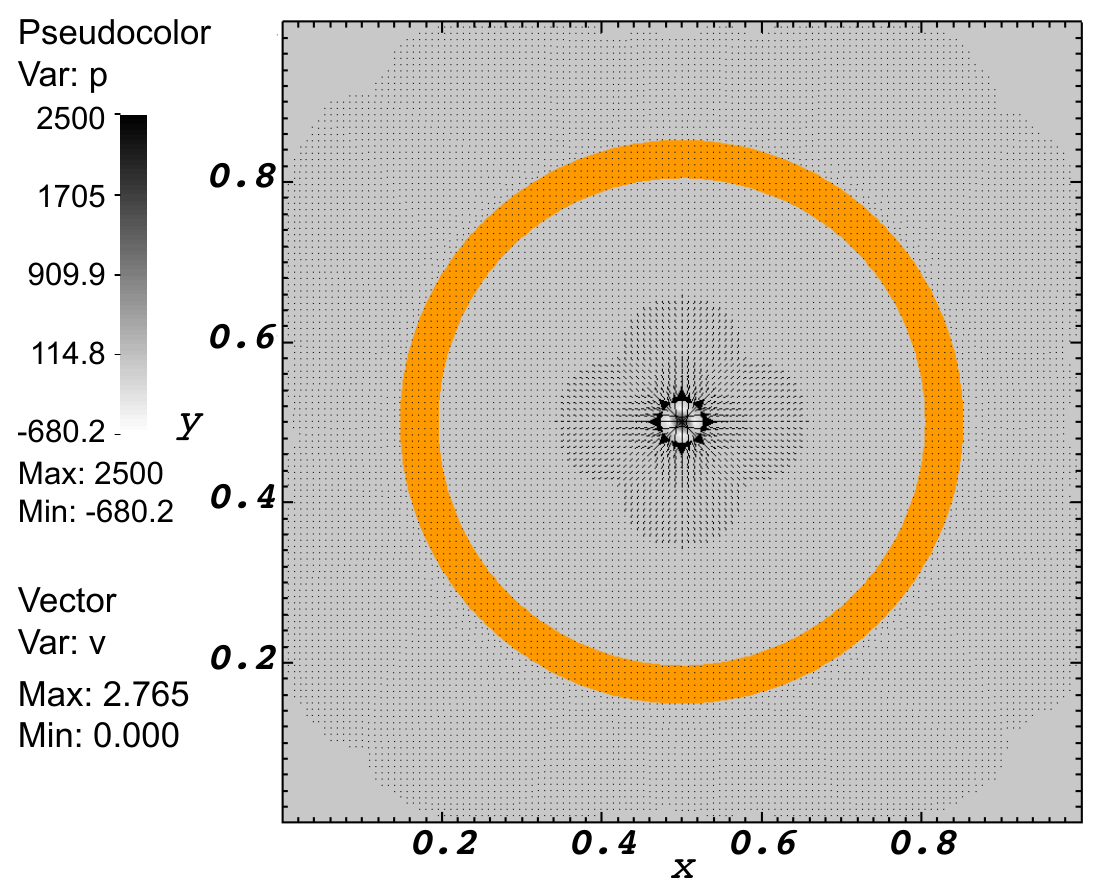}
\label{subfig:pts-comp-solid-t1}
}
\caption{The velocity and the mean normal stress field over the control volume. Also shown is the annulus mesh.}
\end{figure}
As time progresses, the fluid entering the control volume deforms and compresses the annulus, whose configuration for $t=\np[s]{1}$ is shown in Fig.~\ref{subfig:pts-comp-solid-t1}. Referring to Fig.~\ref{fig:pts-comp-solid-error-estimate},
\begin{figure}[htbp]
\centering
\subfigure
{\includegraphics[width=0.48\textwidth, trim = 0.5in 2.6in 0.9in 2.5in, clip=true]{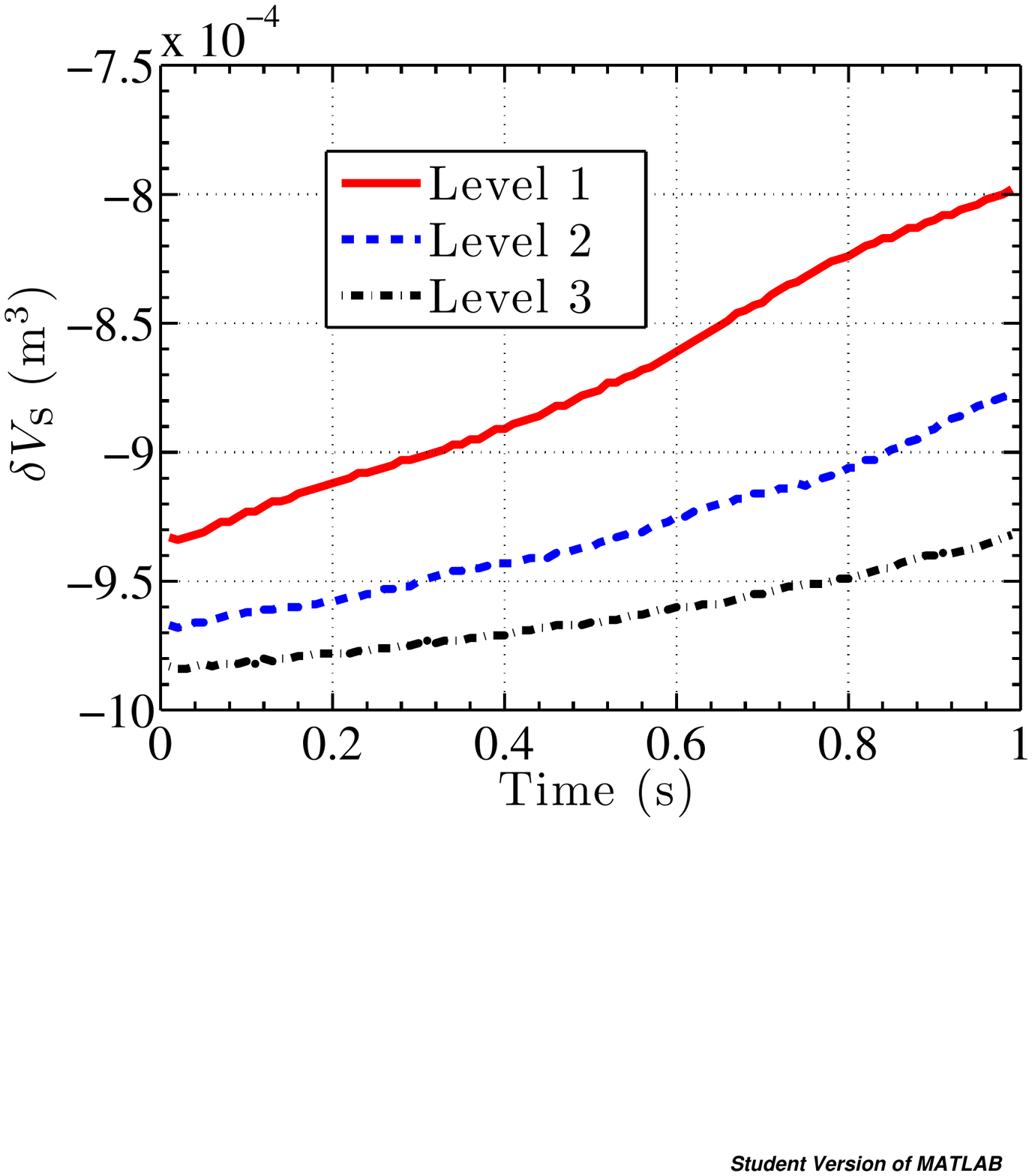}
\label{subfig:pts-comp-solid-error-estimate-1}
}
\subfigure
{\includegraphics[width=0.47\textwidth, trim = 0.6in 2.6in 0.9in 2.5in, clip=true]{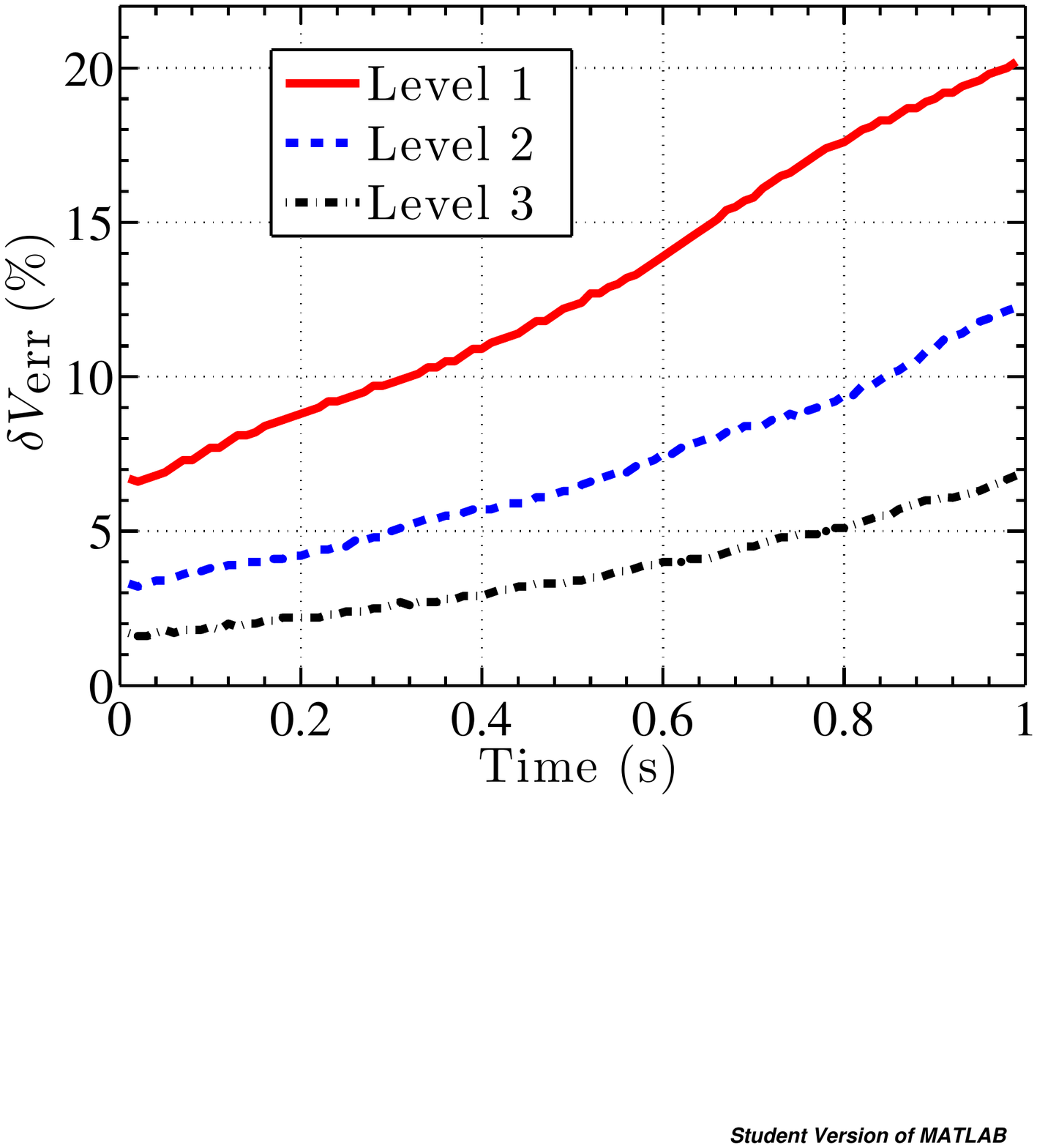}
\label{subfig:pts-comp-solid-error-estimate-2}
}
\caption{The difference between the instantaneous amount of fluid entering due to the source and the change in the area of the annulus. The difference reduces with mesh refinement.}
\label{fig:pts-comp-solid-error-estimate}
\end{figure}
when we look at the difference in the instantaneous amount of fluid entering the control volume and the decrease in the volume of the solid, we see that the difference increases over time. This is not surprising since the mesh of the solid becomes progressively distorted as the fluid emanating from the point source push the inner boundary of the annulus. As expected, the error significantly reduces with the increase in the refinements of the fluid and the solid meshes.

\subsubsection{Elastic bar behind a cylinder}
We now present a second example pertaining to a compressible elastic solid, this time without a viscous component to its behavior.  To test the interaction of a purely elastic compressible object in an incompressible linear viscous flow we consider the two non-steady \acro{FSI} cases discussed by \cite{Turek2006Proposal-for-Nu0}, referred to as \acro{FSI2} and \acro{FSI3}, respectively.  In presenting our results and to facilitate comparisons, we use the same geometry, nondimensionalization, and parameters used by \cite{Turek2006Proposal-for-Nu0}.  Specifically, referring to Fig.~\ref{fig:BarCylinder},
\begin{figure}[htb]
\centering
\includegraphics{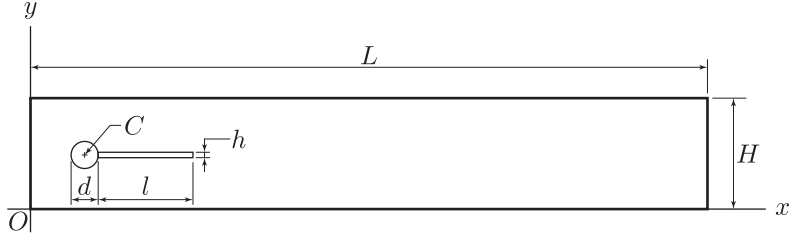}
\caption{Solution domain for the problem of an elastic bar behind a cylinder.}
\label{fig:BarCylinder}
\end{figure}
the system consists of a 2D channel of dimensions $L = \np[m]{2.5}$ and height $H = \np[m]{0.41}$, with a fixed circle $K$ of diameter $d = \np[m]{0.1}$ and centered at $C = (0.2,0.2)\,\text{m}$. The elastic bar attached at the right edge of the circle has length $l = \np[m]{0.35}$ and height $h = \np[m]{0.02}$.

The constitutive response of the bar is that of a de Saint-Venant Kirchhoff material (\citealp{Holzapfel2000Nonlinear-Solid0,Turek2006Proposal-for-Nu0}) so that the viscous component of the stress is equal to zero ($\bv{\sigma}_{\s}^{v} = \tensor{0}$) and the (purely) elastic stress behavior is given by
\begin{equation}
\label{eq: SVK Material}
\bv{\sigma}_{\s}^{e}
=
J^{-1} \tensor{F} \bigl[2 G^{e} \tensor{E} + \lambda^{e} (\trace \tensor{E}) \tensor{I}\bigr] \trans{\tensor{F}}
=
J^{-1} \tensor{F}\biggl[2 G^{e} \tensor{E} + \frac{2 G^{e} \nu^{e}}{1-2 \nu^{e}}
(\trace \tensor{E}) \tensor{I}\biggr] \trans{\tensor{F}},
\end{equation}
where $\tensor{E} = (\trans{\tensor{F}}\tensor{F} - \tensor{I})/2$ is the Lagrangian strain tensor, $G^{e}$ and $\lambda^{e}$ are the Lam{\'e} elastic constants of the immersed solid, and where $\nu^{e} = \lambda^{e}/[2(\lambda^{e} + G^{e})]$ is corresponding Poisson's ratio.

The system is initially at rest.  The boundary conditions are such that there is no slip over the top and bottom surfaces of the channel as well as over the surface of the circle (the immersed solid does not slip relative to the fluid).  At the right end of the channel we impose ``do nothing'' boundary conditions.  Using the coordinate system indicated in Fig.~\ref{fig:BarCylinder}, at the left end of the channel we impose the following distribution of inflow velocity:
\begin{equation}
\label{eq: left BC Turek and Hron}
u_{x} = 1.5 \bar{U} \, 4 y(H - y)/H^{2}
\quad \text{and} \quad
u_{y} = 0,
\end{equation}
where $\bar{U}$ is a constant with dimension of speed.

The constitutive parameters and the the parameter $\bar{U}$ used in the simulations are reported in Table~\ref{tab: Tureck and Hron Parameters}, in which we have also indicated the flow's Reynolds number.
\begin{table}[htbp]\small
\caption{Parameters used in the two non-steady \acro{FSI} cases in \cite{Turek2006Proposal-for-Nu0}.}
\begin{center}
\begin{tabular}{l r r}
\toprule
Parameter & \acro{FSI2} & \acro{FSI3}\\
\midrule
$\rho_{\s}\,(10^{3}\,\text{kg}/\text{m}^{3})$ & $10.0$ & $1.0$ \\
$\nu^{e}$ & 0.4 & 0.4 \\
$G^{e}\,[10^{6}\,\text{kg}/(\text{m}\!\cdot\!\text{s}^{2})]$ & $0.5$ & $2.0$ \\
$\rho_{\f}\,(10^{3}\,\text{kg}/\text{m}^{3})$ & $1.0$ & $1.0$ \\
$\mu_{\f}\,(10^{-3}\,\text{m}^{2}/\text{s})$ & $1.0$ & $1.0$ \\
$\bar{U}\,(\text{m}/\text{s})$ & $1.0$ & $2.0$ \\
$\rho_{\s}/\rho_{\f}$ & $10.0$ & $1.0$ \\
$\reynolds = \bar{U} d/ \mu_{\f}$ & $100.0$ & $200.0$ \\
\bottomrule
\end{tabular}
\end{center}
\label{tab: Tureck and Hron Parameters}
\end{table}

The outcome of the numerical benchmark proposed by \cite{Turek2006Proposal-for-Nu0} is typically expressed (\emph{i}) in terms of the time dependent position of the the midpoint $A$ at the right end of the elastic bar, and (\emph{ii}) in terms of the force acting on the boundary $S$ of the union of the circle $K$ and the elastic bar $B_{t}$ (see Fig.~\ref{fig: THResultsSetUp}).
\begin{figure}[htb]
    \centering
    \includegraphics{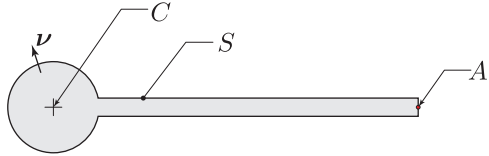}
    \caption{Domain resulting from the union of the elastic bar $B_{t}$ and the fixed circle $K$ with center $C$. $S = \partial(K \cup B_{t})$ denotes the boundary of the domain in question and it is oriented by the unit normal $\bv{\nu}$.  Point $A$ is the midpoint on the right boundary of the elastic bar.}
    \label{fig: THResultsSetUp}
\end{figure}
Denoting by $\ui$ and $\uj$ the orthonormal base vectors associated with the $x$ and $y$ axes, respectively, the force acting on the domain $S$ is 
\begin{equation}
\label{Eq: lift and drag def}
F_{D} \, \ui + F_{L}\,\uj = \int_{S} \bv{\sigma} \bv{\nu} \d{a},
\end{equation}
where $F_{D}$ and $F_{L}$ are the lift and drag, respectively, and where $\bv{\sigma} \bv{\nu}$ is the (time dependent) traction vector acing on $S$.  As remarked by \cite{Turek2006Proposal-for-Nu0} there are several ways to evaluate the right-hand side of Eq.~\eqref{Eq: lift and drag def}.  For example, the traction on the elastic bar could be calculated on the fluid side or on the solid side or even as an average of these values.\footnote{Ideally, the traction value computed on the solid and fluid sides are the same.}  However,  we need to keep in mind that we use a single field for the Lagrange multiplier $p$ and that this field is expected to be discontinuous across the across the (moving) boundary of the immersed object.  In turn, this means that the measure of the hydrodynamic force on $S$ via a direct application of Eq.~\eqref{Eq: lift and drag def} would be adversely affected by the oscillations in the field $p$ across $S$.  With this in mind, referring to the second of Eqs.~\eqref{eq: Balance of mass and momentum}, we observe that a straightforward application of the divergence theorem over the domain $\Omega\setminus(K \cup B_{t})$ yields the following result:
\begin{equation}
\label{Eq: lift and drag def actual}
F_{D} \, \ui + F_{L}\,\uj = \int_{\partial\Omega} \bv{\sigma} \bv{n} \d{a} - \int_{\Omega\setminus(K\cup B_{t})} \rho \biggl\{\bv{b} - \biggl[\frac{\partial\bv{u}}{\partial t} + (\nabla \bv{u}) \bv{u} \biggr] \biggr\} \d{v},
\end{equation}
where $\bv{n}$ denotes the outward unit normal of $\partial\Omega$.  The estimation of the lift and drag over $S$ via Eq.~\eqref{Eq: lift and drag def actual} is significantly less sensitive to the oscillations of the field $p$ near the boundary of the immersed domain and this is the way we have measured $F_{D}$ and $F_{L}$.

For both the \acro{FSI2} and \acro{FSI3} benchmarks, we performed calculations using \np{2992} $\mathcal{Q}_{0}^{2}\vert\mathcal{Q}_{0}^{1}$ elements for the fluid and \np{704} $\mathcal{Q}_{0}^{2}$ elements for the immersed elastic bar.  The number of degrees of freedom distributed over the control volume is \np{24464} for the velocity and \np{3124} for the pressure.  The number of degrees of freedom distributed over the elastic bar is \np{6018}.  The time step size for the \acro{FSI2} benchmark was set to \np[s]{0.005}, whereas the time step size for the \acro{FSI3} results was set to \np[s]{0.001}.

The results for the \acro{FSI2} case are shown in Figures~\ref{fig: FSI 2 displacement BarCylinder} and~\ref{fig: FSI 2 Lift and drag BarCylinder} displaying the components of the displacement of point $A$ and the components of the hydrodynamic force on $S$, respectively.
\begin{figure}[htb]
\centering
\includegraphics{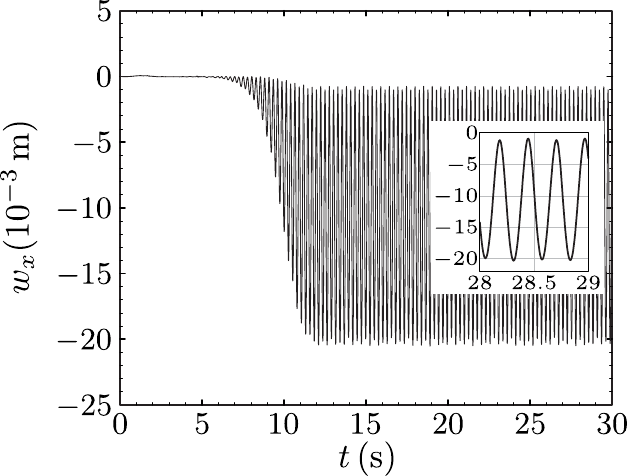}
\quad
\includegraphics{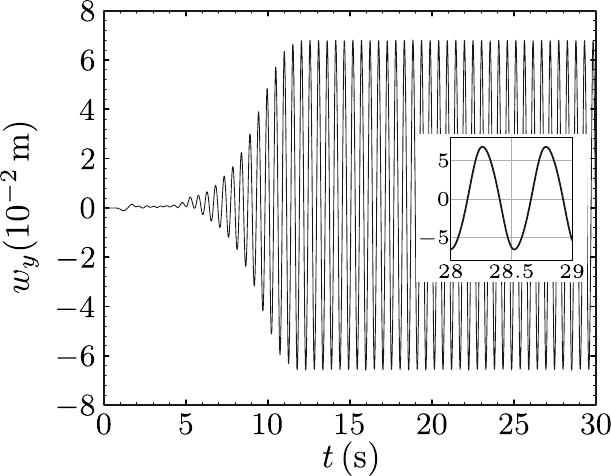}
\caption{Nondimensional displacement of the midpoint at the right end of the elastic bar vs.\ time.  The horizontal and vertical components of the displacement are plotted to the left and to the right, respectively.}
\label{fig: FSI 2 displacement BarCylinder}
\end{figure}
\begin{figure}[htb]
\centering
\includegraphics{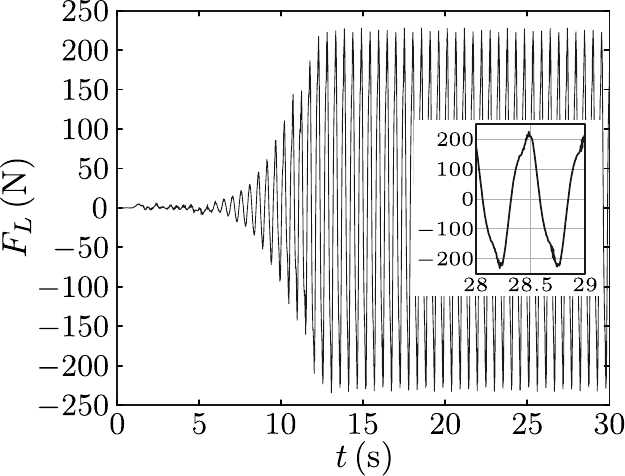}
\quad
\includegraphics{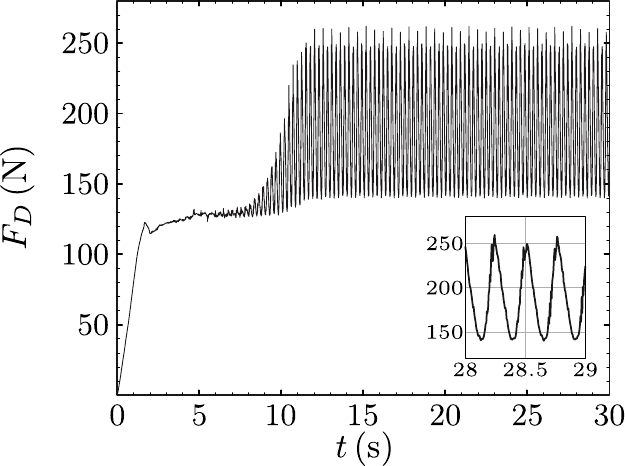}
\caption{Plots of the lift (left) and drag (right) that the fluid exerts on the immersed fixed cylinder and the elastic bar.}
\label{fig: FSI 2 Lift and drag BarCylinder}
\end{figure}
The analogous results for the \acro{FSI3} case are displayed in Figs.~\ref{fig: FSI 3 displacement BarCylinder} and~\ref{fig: FSI 3 Lift and drag BarCylinder}.
\begin{figure}[htb]
\centering
\includegraphics{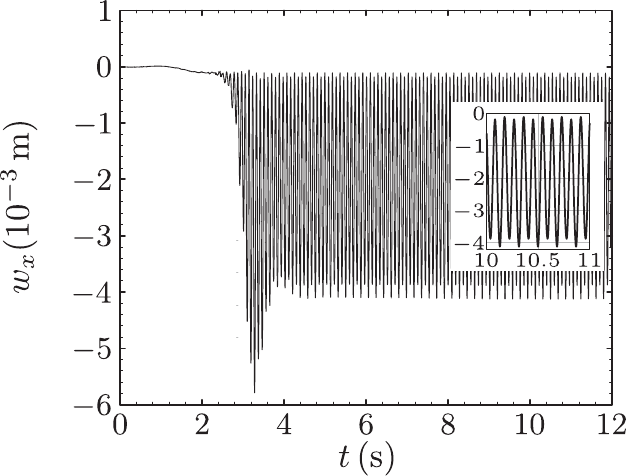}
\quad
\includegraphics{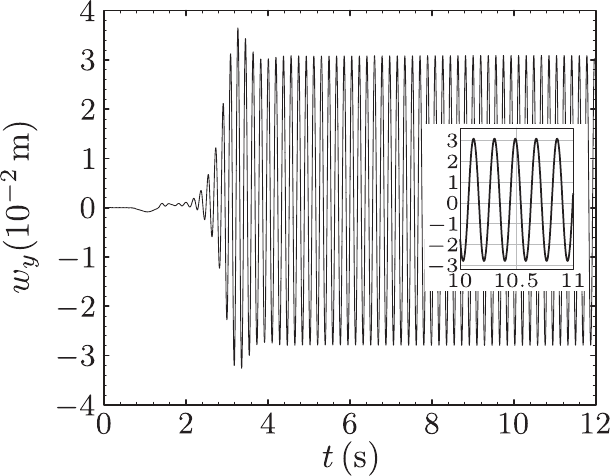}
\caption{Nondimensional displacement of the midpoint at the right end of the elastic bar vs.\ time.  The horizontal and vertical components of the displacement are plotted to the left and to the right, respectively.}
\label{fig: FSI 3 displacement BarCylinder}
\end{figure}
\begin{figure}[htb]
\centering
\includegraphics{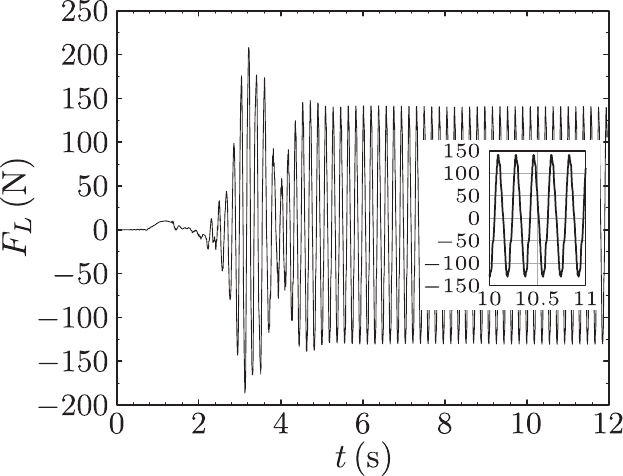}
\quad
\includegraphics{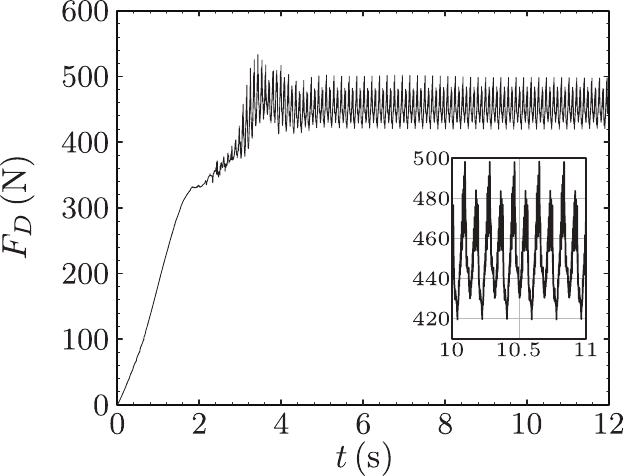}
\caption{Plots of the lift (left) and drag (right) that the fluid exerts on the immersed fixed cylinder and the elastic bar.}
\label{fig: FSI 3 Lift and drag BarCylinder}
\end{figure}
As can be seen in these figures, the displacement values as well as the lift and drag results compare rather favorably with those in the benchmark proposed by \cite{Turek2006Proposal-for-Nu0}, especially when considering that our integration scheme is the implicit Euler method.

\subsection{Flexible 3D bar behind a prismatic rigid obstacle}
Here we present some calculations pertaining to a simple three-dimensional problem.  These results are not intended to reproduce any rigorous benchmarks.  Rather, they provide a snapshot of our current computational capability, which consists of a serial code using \textsc{UMFPACK}~\citep{Davis_2004_UMFPACK} as a direct solver, rather than an optimized parallel code with a carefully designed pre-conditioner. In this sense, the results presented in this section point to the challenges we plan to tackle in the future. A simplified version of the software used to generate the results of this paper is available in~\cite{Heltai2012A-Fully-Coupled0}.

Benchmarks for three-dimensional FSI problems are not as established as those by Turek and Hron for the two-dimensional cases.  Due to the limitations of out current code, we chose the simple three-dimensional problem found in Section~4.4 of the paper by \cite{Wick2011Fluid-Structure-0}.  Referring to Fig.~\ref{fig: 3D benchmark},
\begin{figure}[htb]
\centering
\includegraphics{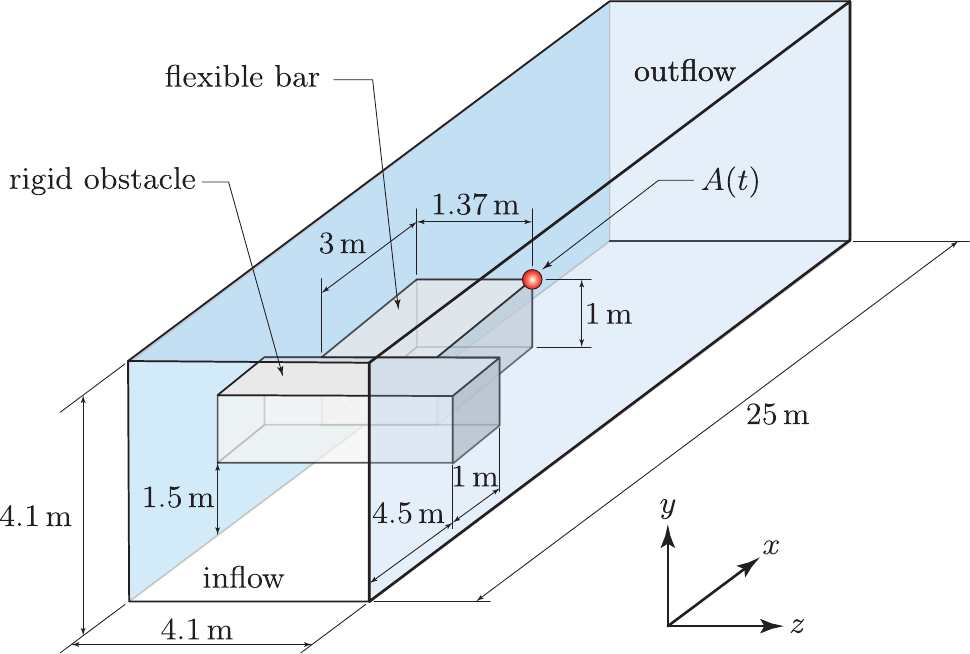}
\caption{Configuration of 3D benchmark.  The results of the simulation consists in tracking the displacement of the point $A(t)$ in the figure and in determining the hydrodynamic forces acting on rigid obstacle and flexible bar system.}
\label{fig: 3D benchmark}
\end{figure}
we consider the motion of an elastic bar attached to a rigid prismatic obstacle across a channel.  Both the channel and the obstacle have square crossections. The obstacle and the elastic bar are not symmetrically located along the $y$ direction within the channel.  In addition, the elastic bar is not positioned symmetrically within the channel along the $z$ direction.  The dimensions of the channel, the rigid obstacle, and the elastic bar attached to the obstacle are shown in the figure.  The stress behavior is purely elastic of the type indicated in Eq.~\eqref{eq: SVK Material}.  The parameters in the simulation are chosen as follows: $\rho_{\f} = \np[kg\!\cdot\!m^{-3}]{1.0}$, $\mu_{\f} = 0.01\,\text{m}^{2}/\text{s}$, $\rho_{\s} = \np[kg\!\cdot\!m^{-3}]{1.0}$, $G^{e} = 500\,\text{kg}/(\text{m}\!\cdot\!\text{s}^{2})$ and $\nu^{e} = 0.4$.  A constant parabolic velocity profile is prescribed at the inlet:
\begin{equation}
\label{eq: 3D example inlet conditions}
u_{x}(0,y,z,t) = 16 U y z (H - y) (H - z) H^{-4},\quad
u_{y}(0,y,z,t) = 0,\quad
u_{z}(0,y,z,t) = 0,
\end{equation}
where $H = \np[m]{4.1}$ and $U = \np[m/s]{0.45}$.  At the outlet ``do nothing'' boundary conditions are imposed. The quantities monitored during the calculations are the components $w_{x}$, $w_{y}$, and $w_{z}$ of the displacement of point $A(t)$, as well as the drag and lift around the obstacle and the elastic bar.  The coordinates of $A$ at the initial time are $(8.5, 2.5, 2.73)\,\text{m}$.

Referring to Table~\ref{tab: 3D discretization}, we have considered two types of discretization, one isotropic and one anisotropic. For each we considered two levels of refinement. Although our code can cope with up to three refinement levels, the solution at each time step for refinement levels higher than two requires several hours of computing time, rendering such grids impractical for unsteady simulations. 
\begin{table}[htbp]\small
\caption{Number of cells and DoFs used in the different simulations involving the deformation of a flexible 3D bar behind a prismatic rigid obstacle.}
\begin{center}
\begin{tabular*} {0.95\textwidth} {@{\extracolsep{\fill}} c c c c c}
\toprule
&
\multicolumn{2}{c}{Solid}&
\multicolumn{2}{c}{Control Volume} 
\\
\cmidrule(r){2-3}\cmidrule(r){4-5}
&
Cells &
DoFs &
Cells &
DoFs
\\
\midrule
Level 1& \np{2} &   \np{132} &  \np{33} &   \np{1434} \\
Level 2& \np{2} &  \np{132} &  \np{150} &  \np{5724} \\
Level 3& \np{16} &  \np{675} & \np{264} & \np{9324} \\
Level 4& \np{16} &  \np{675} & \np{1200} & \np{39432} \\
\bottomrule
\end{tabular*}
\end{center}
\label{tab: 3D discretization}
\end{table}

A summary for all the three-dimensional simulations in terms of the mean values and standard deviation for all tracked quantities is reported in Table~\ref{tab: mean and std analysis for 3D}.
\begin{table}[htbp]\small
\caption{\label{tab: mean and std analysis for 3D}%
Mean and standard deviation values of the components of the displacement of $A(t)$ as well as of the lift and drag over time and as a function of the refinement level.}
\begin{center}
\begin{tabular*} {0.95\textwidth} {@{\extracolsep{\fill}} l c c c c c}
\toprule
&
Level
&
\multicolumn{4}{c}{Time interval size (s)} 
\\
\cmidrule(lr){3-6}
 & & 10 & 20 & 60 & 80
\\
\midrule
\fbox{$w_{x}$: $(\text{mean}, \text{std})\times 10^{3}$\,m} & 1 & 0.269, 21.1 & 0.544, 15.4 & 0.977, 9.10 & 1.38, 8.15
\\
& 2 & 0.370, 8.77 & 0.315, 6.32 & 0.126, 4.37 & ---
\\
& 3 & 1.86, 12.9 & 2.13, 9.14 & --- & ---
\\
& 4 & 0.463, 7.37 & --- & --- & ---
\\
\midrule
\fbox{$w_{y}$: $(\text{mean}, \text{std})\times 10^{3}$\,m} &
1 & 0.00849, 2.08 & 0.865, 2.04 & 2.18, 4.71 & 2.16, 6.16
\\
& 2 & 0.103, 1.07 & 0.0744, 1.89 & 0.333, 4.23 & ---
\\
& 3 & -1.28, 1.34 & -1.60, 1.83 & --- & ---
\\
& 4 & -0.0699, 1.20 & --- & --- & ---
\\
\midrule
\fbox{$w_{z}$: $(\text{mean}, \text{std})\times 10^{3}$\,m} &
1 & -0.857, 2.50 & -1.50, 1.97 & -2.66, 4.70 & -2.24, 5.63
\\
& 2 & 0.0704, 1.26 & -0.176, 3.61 & -0.0857, 7.68 & ---
\\
& 3 & -1.09, 1.51 & -1.26, 1.08 & --- & ---
\\
& 4 & 0.0308, 0.758 & --- & --- & ---
\\
\midrule
\fbox{$F_{L}$: $(\text{mean}, \text{std})$\,N} &
1 & 0.00424, 1.59 & 0.0334, 1.32 & 0.0175, 0.962 & 0.0138, 0.849
\\
& 2 & 0.0264, 1.08 & 0.0317, 0.869 & 0.0407, 0.733 & ---
\\
& 3 & 0.00222, 0.481 & 0.00633, 0.373 & --- & ---
\\
& 4 & 0.0264, 1.08 & --- & --- & ---
\\
\midrule
\fbox{$F_{D}$: $(\text{mean}, \text{std})$\,N} &
1 & 2.65, 82.7 & 2.06, 58.5 & 1.62, 33.8 & 1.58, 29.2
\\
& 2 & 1.73, 47.9 & 1.34, 33.9 & 1.07, 19.6 & ---
\\
& 3 & 1.83, 35.9 & 1.43, 25.4 & --- & ---
\\
& 4 & 1.56, 30.2 & --- & --- & ---
\\
\bottomrule
\end{tabular*}
\end{center}
\end{table}
The results we have obtained with level~1 refinement are shown in Figs.~\ref{fig: 3D w iso 0} and~\ref{fig: 3D CF CL iso 0},
\begin{figure}[htb]
    \centering
    \includegraphics[page=1]{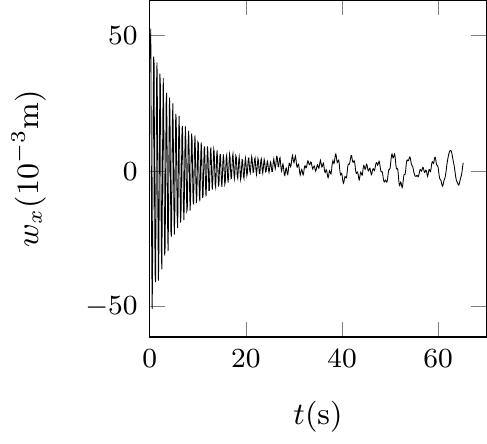}
    \hfill
    \includegraphics[page=2]{Figures/iso_aniso_0.pdf}
    \hfill
    \includegraphics[page=3]{Figures/iso_aniso_0.pdf}
    \caption{Nondimensional displacement of the corner point $A$  at the end of the elastic bar vs.\ time.  The horizontal, vertical, and transversal components of the displacement are plotted from left to right, respectively. Refinement Level 1.}
    \label{fig: 3D w iso 0}
\end{figure}
\begin{figure}[htb]
    \centering
    \includegraphics[page=5]{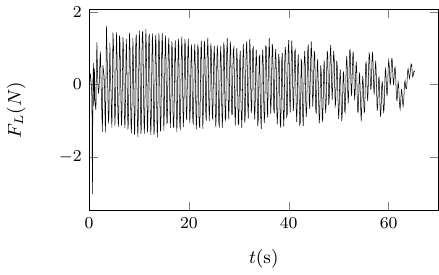}
    \hfill
    \includegraphics[page=4]{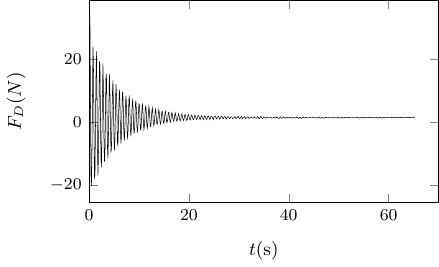}
    \caption{Plots of the lift (left) and drag (right) that the fluid exerts on the primatic section and the elastic bar. Refinement Level 1.}
    \label{fig: 3D CF CL iso 0}
\end{figure}
while those obtained with level~2 refinement are shown in Figs.~\ref{fig: 3D w iso 1} and~\ref{fig: 3D CF CL iso 1}.
\begin{figure}[htb]
    \centering
    \includegraphics[page=6]{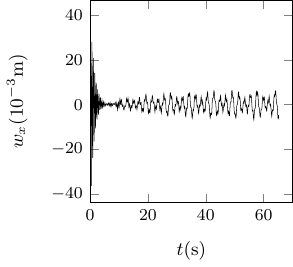}
    \hfill
    \includegraphics[page=7]{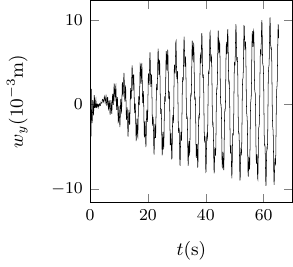}
    \hfill
    \includegraphics[page=8]{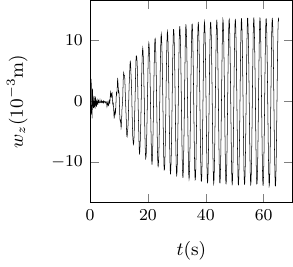}
    \caption{Nondimensional displacement of the corner point $A$  at the end of the elastic bar vs.\ time.  The horizontal, vertical, and transversal components of the displacement are plotted from left to right, respectively. Refinement Level 2.}
    \label{fig: 3D w iso 1}
\end{figure}
\begin{figure}[htb]
    \centering
    \includegraphics[page=10]{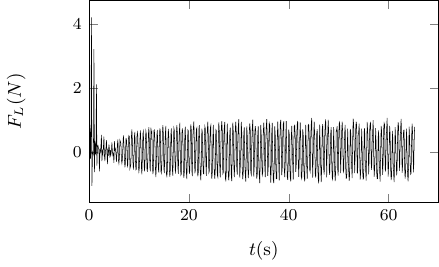}
    \hfill
    \includegraphics[page=9]{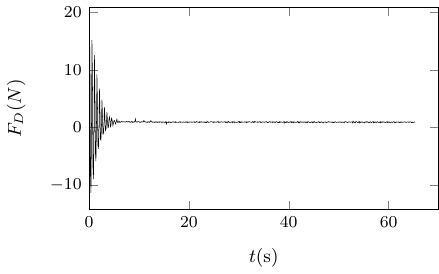}
    \caption{Plots of the lift (left) and drag (right) that the fluid exerts on the primatic section and the elastic bar. Refinement Level 2.}
    \label{fig: 3D CF CL iso 1}
\end{figure}
These simulations display a dynamic response similar to that found in the two-dimensional benchmark problems.  However, no firm conclusion can be drawn other than there is a need for greater refinement to obtain consistent results.  The need for much greater refinement is also evident from the results in \cite{Wick2011Fluid-Structure-0}, where, using a serial code, significantly different outcomes are reported for two consecutive levels of refinement comparable to those used here.  

\section{Summary and Conclusions}
In this paper we have presented the first set of results meant to provide a validation of the fully variational \acro{FEM} approach to an immersed method for \acro{FSI} problems presented in \cite{Heltai2012Variational-Imp0}.  The most important result shows, for the first time, that our proposed immersed method can satisfy the rigorous benchmark tests by \cite{Turek2006Proposal-for-Nu0}.  Our results also show that the proposed approach can be applied to a wide variety of problems in which the immersed solid need not have the same density or the same viscous response as the surrounding fluid.  Furthermore, as shown by the results concerning the lid cavity problem, the proposed computational approach can be applied to problems with very large deformations without any need to adjust the meshes used for either the control volume or the immersed solid.  In the future, we plan to extend these results to include comparative analyses with more established ALE schemes and more extensive and rigorous three-dimensional tests.

\section*{Acknowledgements}
\label{sec:acknowledgements}
The research leading to these results has received specific funding within project OpenViewSHIP, "Sviluppo di un ecosistema computazionale per la progettazione idrodinamica del sistema elica-carena", supported by Regione FVG - PAR FSC 2007-2013, Fondo per lo Sviluppo e la Coesione.

}

\bibliographystyle{FG-AY-bibstyle}               %
\bibliography{rhc_feibm_validation}              %

\end{document}